\documentclass[11pt]{article}
\usepackage{amsmath}
\usepackage{amsfonts}
\usepackage{color}
\usepackage{hyperref}
\newtheorem{theo}{Theorem}[section]
\newtheorem{lem}[theo]{Lemma}
\newtheorem{cor}[theo]{Corollary}

\newtheorem{defi}[theo]{Definition}

\newcommand{\mysection}[1]{\section{#1} \setcounter{equation}{0}}
\newcommand{\proof}{{\sc Proof.} \quad}
\newcommand{\proofc}{{\sc Proof} \ }
\newcommand{\be}{\begin{equation} \label}
\newcommand{\ee}{\end{equation}}
\newcommand{\bea}{\begin{eqnarray}\label}
\newcommand{\eea}{\end{eqnarray}}
\newcommand{\bas}{\begin{eqnarray*}}
\newcommand{\eas}{\end{eqnarray*}}
\newcommand{\bit}{\begin{itemize}}
\newcommand{\eit}{\end{itemize}}
\newcommand{\qed}{\hfill$\Box$ \vskip.2cm}
\newcommand{\nn}{\nonumber}
\newcommand{\R}{\mathbb{R}}
\newcommand{\N}{\mathbb{N}}
\newcommand{\pO}{\partial\Omega}

\newcommand{\eps}{\varepsilon}
\newcommand{\dist}{{\rm dist} \, }
\newcommand{\supp}{{\rm supp} \, }

\newcommand{\wto}{\rightharpoonup}

\newcommand{\hra}{\hookrightarrow}
\newcommand{\io}{\int_\Omega}

\newcommand{\abs}{\\[5pt]}
\newcommand{\ueps}{u_\eps}
\newcommand{\weps}{w_\eps}
\newcommand{\deps}{d_\eps}
\newcommand{\yeps}{y_\eps}
\newcommand{\heps}{h_\eps}
\newcommand{\zeps}{z_\eps}
\newcommand{\eeps}{\eta_\eps}
\newcommand{\deleps}{\delta_\eps}
\newcommand{\tme}{T_{max,\eps}}
\newcommand{\uw}{\underline{w}}
\newcommand{\es}{\eps_\star}
\newcommand{\hatt}{\widehat{T}_\eps}
\newcommand{\set}{{\mathcal{S}}}
\newcommand{\tu}{\widetilde{u}}
\newcommand{\tw}{\widetilde{w}}
\newcommand{\zd}{\zeta_\delta}
\newcommand{\hu}{\widehat{u}}
\newcommand{\hw}{\widehat{w}}
\begin{document}
\enlargethispage{10mm}
\title{Global weak solutions to a strongly degenerate haptotaxis model}
\author{
Michael Winkler\footnote{michael.winkler@math.uni-paderborn.de}\\
{\small Institut f\"ur Mathematik, Universit\"at Paderborn,}\\
{\small 33098 Paderborn, Germany} 
\and
Christina Surulescu\footnote{surulescu@mathematik.uni-kl.de}\\
{\small Technische Universit\"{a}t Kaiserslautern, Felix-Klein-Zentrum f\"{u}r Mathematik,} \\
{\small 67663 Kaiserslautern, Germany}
}

\date{}
\maketitle
\begin{abstract}
\noindent 
We consider a one-dimensional version of a model obtained in \cite{ehs-15} and describing the anisotropic spread of tumor cells in 
a tissue network. The model consists of 
a reaction-diffusion-taxis equation for the density of tumor cells coupled with an ODE for the density of tissue fibers and 
allows for strong degeneracy both in the diffusion and the haptotaxis terms. 
In this setting we prove the global existence of weak solutions to an associated no-flux initial-boundary value problem.\abs
\end{abstract}
%
%
%
%
%
%
%
%
\section{Introduction}\label{intro}
Models with degenerate diffusion in the context of taxis equations have received increased interest during the last decade. They 
describe the dynamics of a cell population in response to a chemoattractant \cite{eewz-14,lw-05,wang-winkler-wrzosek}, moving up the 
gradient of an insoluble signal (haptotaxis) \cite{zhigun_surulescu_uatay}, or performing both chemo- and haptotaxis 
\cite{li-lankeit,tao-winkler-11,zheng-mu-song}. 

In this work we consider a reaction-diffusion-transport-haptotaxis model which is inspired by the effective equations obtained 
in \cite{ehs-15} via parabolic scaling upon starting from a multiscale model for glioma invasion in the anisotropic brain tissue and relying on the setting introduced in  \cite{ehks-15}. 
More precisely, the following PDE-ODE model was considered for the density function $p(t,x,v,y)$ of glioma cells depending on time 
$t$, position $x\in \R^n$, velocity $v\in V:=s\mathbb S^{n-1}$, and density $y\in Y:=(0,R_0)$ of cell surface receptors\footnote{$R_0$ 
denotes the total amount of receptors, assumed to be constant} bound to 
tissue fibers, and for the subcellular dynamics simplified to mass action kinetics of the mentioned receptor binding:
\begin{align}
\partial_tp  +\nabla _x\cdot (vp)+\nabla _y\cdot (G(y,w)p)&=\mathcal L[\lambda ]p+\mathcal P(p)\\
\dot y&=G(y,w).
\end{align}
Thereby, $w(x)$ represents the (macroscopic) volume fraction of tissue, the turning operator  
$\mathcal L[\lambda ]p:=-\lambda (y)p+\int _V\lambda (y)K(x,v,v')p(v')dv'$ 
describes the reorientation of cells due to contact guidance by tissue, and the term 
$\mathcal P(p):=\mu (x,\bar p,v)\int _Y\chi (x,y,y')p(t,x,v,y')w(x)dy'$ models proliferation subsequent to cell-tissue interactions. 
The function $\lambda (y)$ denotes the cell reorientation rate, $K(x,v,v')$ is the turning kernel depending 
on the directional distribution $q(x,v)$ of tissue fibers (obtained from diffusion tensor imaging data), $\mu $ represents the 
proliferation rate depending on the macroscopic cell density $\bar p=\int _V\int _Yp(t,x,v,y)dydv$, and $\chi $ is a kernel 
characterizing the transition from the state $y$ to the state $y'$ during a proliferative action. \abs
An appropriate parabolic scaling led to the macroscopic equation for (an approximation of) the 
tumor cell density:
\begin{equation}\label{eff-macro}
\partial _tu-\nabla \nabla : (\mathbb D_Tu)+\nabla \cdot (a(w)\mathbb D_T\nabla w\ u)=w\mu (x,u)u,
\end{equation}
where $a(w)$ is a function containing both macroscopic and subcellular level information, \\
$\mathbb D_T=\text{const}\int _Vqv\otimes vdv$ is the tumor diffusion tensor encrypting the medical data about the structure of 
brain tissue, and 
\begin{equation}\label{nabla-nabla-eq}
 \nabla \nabla : (\mathbb D_Tu)=\nabla \cdot (\mathbb D_T(x)\nabla u)+\nabla \cdot (\zeta (x)u)
\end{equation}
with the drift velocity $\zeta (x)=\text{const} \int _Vv\otimes v\nabla qdv$. For more details and the precise 
definitions we refer to \cite{ehs-15}. \abs
Equation \eqref{eff-macro} is of the reaction-diffusion-transport-(hapto)taxis type and characterizes the evolution of the tumor 
cell density for a known underlying structure of brain tissue;
in practice, the functions $q$ and $w$ are assessed at a certain time point $t$ from medical data. 
This facilitates both its mathematical analysis and efficient numerical handling, however in fact the 
tumor evolution in a patient also induces dynamical changes in the tissue 
such as e.g.~depletion or remodeling,  which play an essential role in the disease development, see e.g. \cite{bellail,rao} and the references therein. Therefore, a further equation is needed to describe these tissue modifications under the influence of tumor cells. Although in practice it is not feasible from the viewpoint of medical imaging to assess the tissue structure dynamically, by way of model-based predictions relying on such PDE-ODE coupled systems it is possible to use a sequence of just a few 
images in order to obtain via numerical simulations a good approximation of the dynamics over the whole timespan of interest.  \abs
Another issue is related to possible (local) degeneracies of the tumor diffusion tensor $\mathbb D_T(x)$, which 
is particularly relevant e.g.~when modeling resected or irradiated regions of the tumor, where the tissue has been 
depleted as well. 
In the respective domains, this indeed reduces the otherwise
diffusion-dominated PDE \eqref{eff-macro} to a hyperbolic transport equation with nonlinear source term. 
The ma\-the\-ma\-ti\-cal\-ly quite delicate features of such strongly degenerate systems become manifest already in the case when
any taxis or source terms are absent, that is, when $a\equiv 0$ and $\mu\equiv 0$ in (\ref{eff-macro}).
Indeed, in \cite{hillen-painter-winkler} the linear scalar parabolic equation
\begin{equation}\label{eq-hpw}
	\partial _tu=(d_1(y)u)_{xx}+(d_2(y)u)_{yy}, \qquad (x,y)\in \Omega =(0,L_x)\times (0,L_y),\ t>0,
\end{equation}
has been studied, motivated among others by a monoscale model for anisotropic glioma spread in \cite{painter-hillen},
and it was shown there that if
the functions $d_1$ and $d_2$ are smooth and nonnegative and
such that $d_1$ is strictly positive but $d_2$ vanishes precisely in some subinterval $[a,b]$ of $(0,L_y)$,
then solutions to an associated no-flux initial-boundary value problem
asymptotically approach a singular state reflecting concentration of mass
within the degeneracy region $[0,L_x] \times [a,b]$ and extinction outside.\abs
In this paper we intend to provide a first step toward a mathematical understanding of corresponding sytstems when
beyond such strongly degenerate diffusion processes,
further crucial  mechanisms and especially nonlinear haptotaxis are involved.
In order to concentrate on essential aspects of such types of interplay within the framework of a model 
that captures the essential properties but beyond that remains as simple as possible,
we may restrict to the spatially one-dimensional case, in which
the tumor diffusion tensor $\mathbb D_T$ in \eqref{eff-macro} actually reduces to a scalar function.
In the context of a simple evolution law for the haptotactic attractant, particularly neglecting remodeling mechanisms,
this leads to coupled parabolic-ODE systems of the form
\be{v}
	\left\{ \begin{array}{l}
	u_t= \big(d(x)u\big)_{xx} - \big(d(x)u\psi(v) v_x\big)_x, \\[1mm]
	v_t=-uh(v),
	\end{array} \right.
\ee
with given nonnegative functions $d, \psi$ and $h$.\abs
Although in our current 1D setting \eqref{v} the model in \cite{ehs-15} loses most of its 
anisotropy relevance, some of it is retained in the space-dependent 
diffusion and haptotactic sensitivity coefficients. Likewise, the multiscality considered in \cite{ehs-15} and leading to a 
haptotactic coefficient depending on the subcellular dynamics can still be partially retained in this model, in spite of the modified transport term, in which the drift velocity has now a simpler form, yet depending on $d(x)$. The very presence of the 
haptotaxis term is a consequence of taking the receptor binding dynamics into account when describing the evolution of the cell density function on the mesoscopic level and scaling up to the macroscopic one. 
Hence, essential features of the model obtained in \cite{ehs-15} are preserved even in this simplified, 
dimension-reduced setting.\abs
Another related model featuring degenerate diffusion in the context of haptotaxis was proposed and investigated in \cite{zhigun_surulescu_uatay}. The kind of degeneracy considered there is, however, different from the one in this and previous models, as it affects both the diffusion and the haptotaxis coefficients, 
thereby allowing the diffusion to degenerate due to one or both solution components (tumor cell density and tissue density). Unlike the present model, in \cite{zhigun_surulescu_uatay}  there is (apart from the taxis) no other transport term. \abs
{\bf Problem setup and main result.} \quad
In order to make the essential mathematical aspects of (\ref{v}) more transparent, let us 
write (\ref{v}) in a form involving a constant haptotacitc sensitivity, which according to the simple ODE
structure of the second equation therein can readily be achieved on substituting 
$w=\Psi(v)$ with $\Psi(v):=\int_0^v \psi(\sigma)d\sigma$, $v\ge 0$.
Accordingly, in an open bounded interval $\Omega\subset \R$ we will henceforth consider
the initial-boundary value problem
\be{0}
	\left\{ \begin{array}{ll}
	u_t=(d(x)u)_{xx} - (d(x)uw_x)_x + uf(x,u,w),
	\qquad & x\in\Omega, \ t>0, \\[1mm]
	w_t=-ug(w),
	\qquad & x\in\Omega, \ t>0, \\[1mm]
	(d(x)u)_x - d(x)uw_x=0, 
	\qquad & x\in\pO, \ t>0, \\[1mm]
	u(x,0)=u_0(x), \quad w(x,0)=w_0(x),
	\quad & x\in\Omega,
	\end{array} \right.
\ee
with given parameter functions $d:\bar\Omega\to [0,\infty)$, $f:\bar\Omega\times [0,\infty)^2 \to  \R$ and
$g:[0,\infty) \to [0,\infty)$ satisfying
\be{dfg_reg}
	\sqrt{d}\in W^{1,\infty}(\Omega),
	\quad 
	f\in C^1(\bar\Omega\times [0,\infty)^2)
	\quad \mbox{and} \quad
	g\in C^1([0,\infty)),	
\ee
and with prescribed initial data $u_0$ and $w_0$ which are such that
\be{init}
	\left\{ \begin{array}{l}
	0\le u_0\in C^0(\bar\Omega) \mbox{ satisfies } u_0\not\equiv 0 \quad \mbox{and} \\[1mm]
	0\le w_0 \in W^{1,2}(\Omega) \mbox{ has the property that } \io \frac{w_{0x}^2}{g(w_0)} < \infty.
	\end{array} \right.
\ee
As for the parameter functions in (\ref{0}),
throughout our analysis we shall furthermore assume that 
\be{f_upper}
	f(x,u,w) \le \rho(w)
	\quad \mbox{for all $(x,u,w)\in\bar\Omega\times [0,\infty)^2$ \quad with some nondecreasing } 
	\rho: [0,\infty) \to [0,\infty),
\ee
and that there exists $\delta>0$ such that writing
\be{M}
	M:=\|w_0\|_{L^\infty(\Omega)}+\delta,
\ee
we have
\be{g}
	g(0)=0, 
	\qquad
	g(w)>0 \quad \mbox{for all } w\in (0,M]
	\qquad \mbox{and} \qquad
	g'(w) > 0 \quad \mbox{for all } w\in [0,M]
\ee
as well as 
\be{g0}
	\liminf_{w\searrow 0} \frac{g'(w)}{g(w)}>0,
\ee
whence in particular there exist $\Gamma>0$ and $\gamma>0$ fulfilling
\be{g_upper}
	g(w) \le \Gamma w
	\quad \mbox{for all } w\in [0,M]
\ee
 and
\be{gamma}
	\frac{g'(w)}{g(w)} \ge \gamma 
	\qquad \mbox{for all } w\in (0,M].
\ee
Beyond the analytically simplest case obtained on letting
\bas
	g(w)=w, \qquad w\ge 0,
\eas
this inter alia includes more general choices such as 
\bas
	g(w)=w(1-w),	\qquad w\ge 0,
\eas
upon which via the substitution $w=\frac{v}{1+v}$, on the set of solutions fulfilling $v<1$ the system
(\ref{0}) becomes formally equivalent to a corresponding initial-boundary value problem for the special version
\be{vs}
	\left\{ \begin{array}{l}
	u_t= \big(d(x)u\big)_{xx} - \Big(\frac{d(x)u}{(1+v)^2}v_x\Big)_x, \\[1mm]
	v_t=-uv,
	\end{array} \right.
\ee
of (\ref{v}), as proposed in \cite{zhigun_surulescu_uatay} for modeling tumor invasion in a tissue network, thereby paying increased attention to the form of the haptotaxis coefficient. Specifically, the latter accounts for microscopic cell-tissue interactions, which --besides having a haptotaxis term at all-- retains a supplementary trace of multiscality in our macroscopic model, although in a rather indirect way, as we do not explicitly couple some ODE for receptor binding kinetics to the two PDEs for $u$ and $v$. The presence of $d(x)$ in both diffusion/transport and haptotaxis coefficients is motivated by the deduction in \cite{ehs-15}.\abs 
The main results of our analysis indicate that even in this general setting, thus allowing for virtually arbitrary
strength of degeneracies in diffusion, haptotactic cross-diffusion does not result in a finite-time collapse of
solutions into e.g.~persistent Dirac-type singularities. 
More precisely, let us introduce the following solution concept to pursued below,
in which we use the abbreviation $\{d>0\} := \big\{ x\in \bar\Omega \ \Big| \ d(x)>0 \big\}$
which along with a corresponding definition of $\{d=0\}$ will frequently be used throughout the sequel.
\begin{defi}\label{defi_weak}
  A pair $(u,w)$ of nonnegative functions
  \be{w1}
	\left\{ \begin{array}{l}
	u \in L^1_{loc}(\bar\Omega\times [0,\infty)), \\[1mm]
	w\in L^\infty_{loc}(\bar\Omega\times [0,\infty)) \cap L^1_{loc}([0,\infty);W^{1,1}(\{d>0\}))
	\end{array} \right.
  \ee
  satisfying
  \be{w11}
	uf(\cdot,u,w) \in L^1_{loc}(\bar\Omega\times [0,\infty))
	\qquad \mbox{and } \qquad
	ug(w) \in L^1_{loc}(\bar\Omega\times [0,\infty))
  \ee
  as well as
  \be{w2}
	duw_x \in L^1_{loc}([0,\infty);L^1(\{d>0\}))
  \ee
  will be called a {\em global weak solution} of (\ref{0}) if
  \bea{w3}
	- \int_0^\infty \io u\varphi_t - \io u_0 \varphi(\cdot,0)
	= \int_0^\infty \int_{\{d>0\}} du\varphi_{xx}
	+ \int_0^\infty \int_{\{d>0\}} duw_x \varphi_x + \int_0^\infty \io uf(\cdot,u,w) \varphi
  \eea
  for all $\varphi\in C_0^\infty(\bar\Omega\times [0,\infty))$ such that $\varphi_x=0$ on $\pO \times (0,\infty)$
  and
  \be{w4}
	\int_0^\infty \io w\varphi_t + \io w_0\varphi(\cdot,0)
	= \int_0^\infty \io ug(w) \varphi
  \ee	
  for all $\varphi\in C_0^\infty(\Omega\times [0,\infty))$.
\end{defi}
Within this framework, a global solution of (\ref{0}) can always be constructed:
\begin{theo}\label{theo33}
  Suppose that $\Omega\subset \R$ is a bounded interval, and that $u_0, w_0$, $d, f$ and $g$ satisfy
  (\ref{init}), (\ref{dfg_reg}), (\ref{f_upper}) and (\ref{g}).
  Then (\ref{0}) possesses at least one global weak solution in the sense specified in Definition \ref{defi_weak} below.
\end{theo}
This paper is organized as follows: In Section \ref{regularized-problem} we introduce a regularized version of the degenerate
problem, for which some useful properties are obtained. Section \ref{entropy} is concerned with studying an entropy functional which allows to deduce a quasi-dissipative property of the regularized system, inter alia asserting 
global existence of its solution. 
Some precompactness and regularity properties of terms involved in that system follow in Sections \ref{precomp} 
and \ref{reg-prop}, respectively, succeeded in Section \ref{reg-time} by regularity features of corresponding 
time derivatives. 
Sections \ref{limits-d-positive} and \ref{further-limits} provide convergence properties of the approximate solution in the region with no degeneracy; further properties of the respective limits are obtained in Section \ref{sol-prop}. 
Finally, Section \ref{the-proof} concludes the existence proof for the strongly degenerate problem (\ref{0}).
\mysection{Regularized problems and their basic properties}\label{regularized-problem}
In order to prepare the construction of an appropriate family of non-degenerate approximations of (\ref{0}),
according to the nonnegativity of $d$ and the inclusion $\sqrt{d}\in W^{1,\infty}(\Omega)$ 
we may first choose $(\deps)_{\eps\in (0,1)} \subset C^3(\bar\Omega)$ in such a way that 
$d_{\eps x}=0$ on $\pO$ and that 
with some $K_1>0$, for each $\eps\in (0,1)$ we have
\be{d_pos}
	\sqrt{\eps} \le \deps(x) \le \|d\|_{L^\infty(\Omega)}+1
	\qquad \mbox{for all } x\in\bar\Omega,
\ee
as well as	
\be{K1}
	\frac{d_{\eps x}^2(x)}{\deps(x)} \le K_1
	\qquad \mbox{for all } x\in\bar\Omega,
\ee
and such that moreover
\be{d_conv}
	\deps\to d
	\quad \mbox{in } L^\infty(\Omega)
	\qquad \mbox{as } \eps\searrow 0
\ee
and
\be{dx_conv}
	d_{\eps x} \to d_x
	\quad \mbox{a.e.~in } \Omega
	\qquad \mbox{as } \eps\searrow 0.
\ee
We next note that according to (\ref{g}) it is possible to fix $\eps_0\in (0,1)$ such that $g(M)>\eps_0$,
whereupon with $\delta$ as introduced in the course of the definition (\ref{M}) of $M$,
for each $\eps\in (0,\eps_0)$ we can choose $\deleps \in (0,\delta^2)$ such that
\be{g_pos}
	g(w) \ge \eps
	\quad \mbox{for all } w\in [\deleps,M],
\ee
and such that moreover $\deleps\to 0$ as $\eps\searrow 0$.
It is then easy to see that one can find $(\eeps)_{\eps\in (0,\eps_0)} \subset (0,1)$ with the two properties that
\be{eta_div}
	\eeps \ln \frac{1}{\sqrt{\deleps}} \to + \infty
	\qquad \mbox{as } \eps\searrow 0,
\ee
and that
\be{eta_conv}
	\eeps\to 0
	\qquad \mbox{as } \eps\searrow 0;
\ee
indeed, it can readily be checked that this can be achieved on choosing
\bas
	\eeps:=\frac{\ln \ln \frac{A}{\sqrt{\deleps}}}{\ln \frac{A}{\sqrt{\deleps}}},
	\qquad \eps\in (0,\eps_0),
\eas
with some suitably large $A>0$.
For $\eps\in (0,\eps_0)$, we then let 
\be{w0eps}
	w_{0\eps}(x):=w_0(x)+\sqrt{\deleps},
	\qquad x\in\bar\Omega,
\ee
and consider the regularized variant of (\ref{0}) given by
\be{0eps}
	\left\{ \begin{array}{ll}
	u_{\eps t} = (\deps\ueps)_{xx} - \Big(\deps\frac{\ueps}{(1+\eeps\ueps)^2} w_{\eps x}\Big)_x 
	+ \ueps f(x,\ueps,\weps),
	\quad & x\in\Omega, \ t>0, \\[1mm]
	w_{\eps t} = \eps \Big(\frac{w_{\eps x}}{\sqrt{g(\weps)}} \Big)_x - \frac{\ueps}{1+\eeps\ueps} g(\weps),
	\quad & x\in\Omega, \ t>0, \\[1mm]
	u_{\eps x}=w_{\eps x}=0,
	\quad & x\in\pO, \ t>0, \\[1mm]
	\ueps(x,0)=u_0(x), \qquad \weps(x,0)=w_{0\eps}(x),
	\quad & x\in\Omega,
	\end{array} \right.
\ee
Due to the additionally introduced artificial diffusion in the equation for $\weps$ each of these problems
can be viewed as a variant of the well-studied Keller-Segel chemotaxis system;
in fact, as can be seen by straightforward adaptation of arguments well-established in the analysis of chemotaxis problems
(\cite{amann}, \cite{horstmann_win}, \cite{win_CPDE2010}), all
these problems allow for local-in-time classical solutions which enjoy a favorable extensibility criterion:
\begin{lem}\label{lem_loc}
  For each $\eps\in (0,\eps_0)$, there exist $\tme\in (0,\infty]$ and nonnegative functions 
  \bas
	\left\{ \begin{array}{l}
	\ueps \in C^0(\bar\Omega\times [0,\tme)) \cap C^{2,1}(\bar\Omega\times (0,\tme)), \\[1mm]
	\weps \in C^0([0,\tme);W^{1,2}(\Omega)) \cap C^{2,1}(\bar\Omega\times (0,\tme)),
	\end{array} \right.
  \eas
  which solve (\ref{0eps}) in the classical sense in $\Omega\times (0,\tme)$, and which are such that
  \be{ext}
	\mbox{if $\tme<\infty$, then }
	\limsup_{t\nearrow \tme} \bigg( \|\ueps(\cdot,t)\|_{L^\infty(\Omega)} + \|\weps(\cdot,t)\|_{W^{1,2}(\Omega)} 
	+ \Big\|\frac{1}{g(\weps(\cdot,t))}\Big\|_{L^\infty(\Omega)} \bigg)
	= \infty.
  \ee
\end{lem}
Let us first collect some basic properties of these solutions. We first assert some useful pointwise upper and lower
bounds for $\weps$.
\begin{lem}\label{lem1}
  Let $\eps\in (0,\eps_0)$. Then
  \be{1.1}
	\weps(x,t) \le M
	\qquad \mbox{for all $x\in\Omega$ and } t\in (0,\tme)
  \ee
  and
  \be{1.2}
	\weps(x,t) \ge \sqrt{\deleps} e^{-\frac{\Gamma}{\eta _\eps} t}
	\qquad \mbox{for all $x\in\Omega$ and } t\in (0,\tme),
  \ee
  where $\Gamma>0$ is as in (\ref{g_upper}).
\end{lem}
\proof
  Since according to our choices of $\deleps$ and $M$ we have
  \bas
	\weps(x,0)=w_0(x)+\sqrt{\deleps} \le \|w_0\|_{L^\infty(\Omega)}+ \delta=M
	\qquad \mbox{for all } x\in\Omega,
  \eas
  the inequality in (\ref{1.1}) immediately results from the maximum principle applied to the second equation 
  in (\ref{0eps}).
  As a consequence thereof, in view of (\ref{g}) we know that $g'(\weps)\ge 0$ in $\Omega\times (0,\tme)$, whence
  \bas
	\frac{\ueps}{1+\eeps\ueps} g(\weps) \le \frac{1}{\eeps} g(\weps)
	\qquad \mbox{in } \Omega\times (0,\tme),
  \eas
  so that using (\ref{0eps}) and (\ref{g_upper}) we see that
  \bas
	w_{\eps t} \ge \Big(\frac{w_{\eps x}}{\sqrt{g(\weps)}}\Big)_x - \frac{\Gamma}{\eeps} \weps
	\qquad \mbox{in } \Omega\times (0,\tme).
  \eas
  Since
  \bas
	\uw(x,t):=\sqrt{\deleps} e^{-\frac{\Gamma}{\eeps} t},
	\qquad x\in\bar\Omega, \ t\ge 0,
  \eas
  satisfies
  \bas
	\uw_t - \Big(\frac{\uw_x}{\sqrt{g(\weps)}}\Big)_x + \frac{\Gamma}{\eeps} \uw =0
	\qquad \mbox{in } \Omega\times (0,\infty)
  \eas
  and $\frac{\partial \uw}{\partial\nu}=0$ on $\pO\times (0,\infty)$ as well as
  \bas
	\uw(x,0)=\sqrt{\deleps} \le \weps(x,0)
	\qquad \mbox{for all } x\in\Omega
  \eas
  by (\ref{w0eps}), the comparison principle therefeore ensures that $\weps \ge \uw$ in $\Omega\times (0,\tme)$
  and that thus also (\ref{1.2}) is valid.
\qed
Using the latter along with (\ref{f_upper}), we easily obtain the following information on the evolution
of $\io \ueps$.
\begin{lem}\label{lem2}
  With $M$ as defined in (\ref{M}) and $\rho$ taken from (\ref{f_upper}), we have
  \be{2.1}
	\frac{d}{dt} \io \ueps
	\le \rho(M) \io \ueps
	\qquad \mbox{for all } t\in (0,\tme)
  \ee
  and
  \be{2.2}
	\io \ueps(\cdot,t) \le \bigg\{ \io u_0 \bigg\} \cdot e^{\rho(M)t}
	\qquad \mbox{for all } t\in (0,\tme).
  \ee
\end{lem}
\proof
  We integrate the first equation in (\ref{0eps}) and use (\ref{f_upper}) together with (\ref{1.1}) to find that
  \bas
	\frac{d}{dt} \io \ueps 
	= \io \ueps f(x,\ueps,\weps)
	\le \io \ueps \rho(\weps)
	\le \rho(M) \io \ueps
	\qquad \mbox{for all } t\in (0,\tme),
  \eas
  and that hence (\ref{2.1}) holds, from which in turn (\ref{2.2}) results upon integration in time.
\qed
\mysection{Implications of an entropy-like structure}\label{entropy}
Now the core of our approach consists in the detection of a favorable quasi-dissipative property of the system (\ref{0eps})
which can be revealed by following the well-established strategy of considering the time evolution of
a functional that combines a logarithmic entropy of the cell distribution with a properly chosen summand annihilating
the correspondingly obtained cross-diffusive interaction integral.
In order to clarify which precise form the latter takes in the context of the approximate problems (\ref{0eps}),
let us begin by separately tracking the logarithmic entropy.
\begin{lem}\label{lem3}
  Let $\rho, M$ and $K_1$ be as introduced in (\ref{f_upper}), (\ref{M}) and (\ref{K1}).
  Then for each $\eps\in (0,\eps_0)$ we have
  \bea{3.1}
	\frac{d}{dt} \io \ueps \ln \ueps
	+ \frac{1}{2} \io \deps \frac{u_{\eps x}^2}{\ueps}
	&\le& \io \deps \frac{u_{\eps x}}{(1+\eeps\ueps)^2} w_{\eps x} 
	+ \Big(\rho(M)+\frac{K_1}{2}\Big) \cdot \bigg\{ \io u_0 \bigg\} \cdot e^{\rho(M)t} \nn\\
	& & + \io \ueps \ln \ueps \cdot f(x,\ueps,\weps) 
	\qquad \mbox{for all } t\in (0,\tme).
  \eea
\end{lem}
\proof
  Since $\ueps$ is positive in $\bar\Omega\times (0,\tme)$ by the strong maximum principle, we may multiply the first 
  equation in (\ref{0eps}) by $\ln \ueps$ and integrate by parts to see that
  \bea{3.2}
	\frac{d}{dt} \io \ueps\ln \ueps
	&=& \io u_{\eps t} \ln \ueps + \frac{d}{dt} \io \ueps \nn\\
	&=& - \io (\deps \ueps)_x \cdot \frac{u_{\eps x}}{\ueps}
	+ \io \deps \frac{u_{\eps x}}{(1+\eeps\ueps)^2} w_{\eps x} 
	+ \io \ueps \ln \ueps \cdot f(x,\ueps,\weps) \nn\\
	& & + \frac{d}{dt} \io \ueps
	\qquad \mbox{for all } t\in (0,\tme).
  \eea
  Here by Lemma \ref{lem2} we have
  \be{3.3}
	\frac{d}{dt} \io \ueps \le \rho(M) \io \ueps
	\qquad \mbox{for all } t\in (0,\tme),
  \ee
  and using Young's inequality and (\ref{K1}) we obtain
  \bea{3.4}
	- \io (\deps \ueps)_x \cdot \frac{u_{\eps x}}{\ueps} 
	&=& - \io \deps \frac{u_{\eps x}^2}{\ueps}
	- d_{\eps x} u_{\eps x} \nn\\
	&\le& - \frac{1}{2} \io \deps \frac{u_{\eps x}^2}{\ueps} + \frac{1}{2} \io \frac{d_{\eps x}^2}{\deps} \ueps \nn\\
	&\le& - \frac{1}{2} \io \deps \frac{u_{\eps x}^2}{\ueps} + \frac{K_1}{2} \io \ueps
	\qquad \mbox{for all } t\in (0,\tme).
  \eea
  Since 
  \bas
	\Big(\rho(M) + \frac{K_1}{2}\Big) \io \ueps
	\le \Big(\rho(M)+\frac{K_1}{2}\Big)\cdot\bigg\{ \io u_0\bigg\} \cdot e^{\rho(M)t}
	\qquad \mbox{for all } t\in (0,\tme)
  \eas
  by Lemma \ref{lem2}, combining (\ref{3.2})-(\ref{3.4}) thus yields (\ref{3.1}).
\qed
Thanks to a favorable exact relationship between the approximate signal absorption rate 
$0\le u\mapsto \frac{u}{1+\eeps u}$ and the tactic sensitivity 
$0\le \frac{1}{(1+\eeps u)^2}$ in (\ref{0eps}),
an exact compensation of the first summand on the right of (\ref{3.1}) can be achieved on complementing the above
by the following.	
\begin{lem}\label{lem4}
  With $K_1$ as in (\ref{K1}), we have
  \bea{4.1}
	& & \hspace*{-30mm}
	\frac{1}{2} \frac{d}{dt} \io \deps \frac{w_{\eps x}^2}{g(\weps)}
	+ \frac{\eps}{2} \io \deps \cdot \frac{1}{\sqrt{g(\weps)}} \cdot \Big(\frac{w_{\eps x}}{\sqrt{g(\weps)}} \Big)_{x}^2
	+ \frac{1}{2} \io \deps \cdot \frac{\ueps}{1+\eeps\ueps} \cdot \frac{g'(\weps)}{g(\weps)} w_{\eps x}^2 \nn\\
	&\le& - \io \deps \frac{u_{\eps x}}{(1+\eeps\ueps)^2} w_{\eps x}
	+ \frac{\eps K_1}{2} \io \frac{w_{\eps x}^2}{\sqrt{g(\weps)}^3} 
	\qquad \mbox{for all } t\in (0,\tme)
  \eea
  whenever $\eps\in (0,\eps_0)$.
\end{lem}
\proof
  Using that $g(\weps)$ is positive in $\bar\Omega\times [0,\tme)$ due to Lemma \ref{lem1} and (\ref{g}), on the basis
  of the second equation in (\ref{0eps}) we compute
  \bea{4.2}
	\frac{d}{dt} \io \deps \frac{w_{\eps x}^2}{g(\weps)}
	&=& 2 \io \deps \frac{1}{g(\weps)} w_{\eps x} w_{\eps x t}
	- \io \deps \frac{g'(\weps)}{g^2(\weps)} w_{\eps x}^2 w_{\eps t} \nn\\
	&=& 2\io \deps \frac{1}{g(\weps)} w_{\eps x} \cdot \bigg\{ \eps \Big(\frac{w_{\eps x}}{\sqrt{g(\weps)}} \Big)_{xx}
	- \Big(\frac{\ueps}{1+\eeps\ueps} g(\weps)\Big)_x \bigg\} \nn\\
	& & - \io \deps \frac{g'(\weps)}{g^2(\weps)} w_{\eps x}^2 \cdot \bigg\{
	\eps \Big(\frac{w_{\eps x}}{\sqrt{g(\weps)}}\Big)_x 
	- \frac{\ueps}{1+\eeps\ueps} g(\weps) \bigg\} \nn\\
	&=& -2\eps \io \Big(\deps \frac{1}{g(\weps)} w_{\eps x}\Big)_x \cdot
	\Big(\frac{w_{\eps x}}{\sqrt{g(\weps)}} \Big)_x \nn\\
	& & - 2\io \deps \frac{1}{g(\weps)} w_{\eps x} \cdot \Big(\frac{\ueps}{1+\eeps\ueps} g(\weps)\Big)_x \nn\\
	& & - \eps \io \deps \frac{g'(\weps)}{g^2(\weps)} w_{\eps x}^2 \cdot \Big(\frac{w_{\eps x}}{\sqrt{g(\weps)}}\Big)_x
		\nn\\
	& & + \io \deps \frac{\ueps}{1+\eeps\ueps} \frac{g'(\weps)}{g(\weps)} w_{\eps x}^2
	\qquad \mbox{for all } t\in (0,\tme).
  \eea
  Here we expand 
  \bas
	- 2\io \deps \frac{1}{g(\weps)} w_{\eps x} \cdot \Big(\frac{\ueps}{1+\eeps\ueps} g(\weps)\Big)_x 
	&=& -2 \io \deps \frac{u_{\eps x}}{(1+\eeps\ueps)^2} w_{\eps x} \nn\\
	& & - 2\io \deps \frac{\ueps}{1+\eeps\ueps} \frac{g'(\weps)}{g(\weps)} w_{\eps x}^2
	\qquad \mbox{for all } t\in (0,\tme),
  \eas
  so that 
  \bea{4.3}
	& & \hspace*{-20mm}
	- 2\io \deps \frac{1}{g(\weps)} w_{\eps x} \cdot \Big(\frac{\ueps}{1+\eeps\ueps} g(\weps)\Big)_x 
	+ \io \deps \frac{\ueps}{1+\eeps\ueps} \frac{g'(\weps)}{g(\weps)} w_{\eps x}^2 \nn\\
	&=& - 2\io \deps \frac{u_{\eps x}}{(1+\eeps\ueps)^2} w_{\eps x}
	- \io \deps \frac{\ueps}{1+\eeps\ueps} \frac{g'(\weps)}{g(\weps)} w_{\eps x}^2
	\qquad \mbox{for all } t\in (0,\tme).
  \eea
  We next use the identity   
  \be{4.4}
	w_{\eps xx} = \sqrt{g(\weps)} \cdot \Big(\frac{w_{\eps x}}{\sqrt{g(\weps)}}\Big)_x
	+ \frac{g'(\weps)}{2g(\weps)} w_{\eps x}^2
	\qquad \mbox{in } \Omega\times (0,\tme)
  \ee
  to rewrite
  \bas
	\Big(\deps \frac{1}{g(\weps)} w_{\eps x}\Big)_x
	&=& \deps \frac{1}{g(\weps)} w_{\eps xx} 
	-\deps \frac{g'(\weps)}{g^2(\weps)} w_{\eps x}^2
	+ d_{\eps x} \frac{1}{g(\weps)} w_{\eps x} \\
	&=& \deps \frac{1}{\sqrt{g(\weps)}} \Big(\frac{w_{\eps x}}{\sqrt{g(\weps)}}\Big)_x
	- \frac{1}{2} \deps \frac{g'(\weps)}{g^2(\weps)} w_{\eps x}^2
	+ d_{\eps x} \frac{1}{g(\weps)} w_{\eps x}
	\qquad \mbox{in } \Omega\times (0,\tme),
  \eas
  so that on the right-hand side of (\ref{4.2}) we can employ Young's inequality to see that
  \bas
	& & \hspace*{-20mm}
	-2\eps \io \Big(\deps \frac{1}{g(\weps)} w_{\eps x}\Big)_x \cdot
	\Big(\frac{w_{\eps x}}{\sqrt{g(\weps)}} \Big)_x 
	- \eps \io \deps \frac{g'(\weps)}{g^2(\weps)} w_{\eps x}^2 \cdot \Big(\frac{w_{\eps x}}{\sqrt{g(\weps)}}\Big)_x \nn\\
	&=& -2\eps \io \deps \frac{1}{\sqrt{g(\weps)}} \Big( \frac{w_{\eps x}}{\sqrt{g(\weps)}}\Big)_x^2
	- 2\eps \io d_{\eps x} \frac{1}{g(\weps)} w_{\eps x} \cdot \Big(\frac{w_{\eps x}}{\sqrt{g(\weps)}}\Big)_x \\
	&\le& -\eps \io \deps \frac{1}{\sqrt{g(\weps)}} \Big( \frac{w_{\eps x}}{\sqrt{g(\weps)}}\Big)_x^2
	+ \eps\io \frac{d_{\eps x}^2}{\deps} \frac{1}{\sqrt{g(\weps)}} w_{\eps x}^2 \\
	&\le& -\eps \io \deps \frac{1}{\sqrt{g(\weps)}} \Big( \frac{w_{\eps x}}{\sqrt{g(\weps)}}\Big)_x^2
	+ \eps K_1 \io \frac{1}{\sqrt{g(\weps)}} w_{\eps x}^2
	\qquad \mbox{for all } t\in (0,\tme),
  \eas
  again due to (\ref{K1}).
  In conjunction with (\ref{4.3}) and (\ref{4.2}) this yields (\ref{4.1}).
\qed
In fact, combining the latter two lemmata yields a quasi-entropy inequality, the essential implications
of which can be summarized as follows.
\begin{lem}\label{lem5}
  Let $T>0$. Then there exist $\es(T)\in (0,\eps_0)$ and $C(T)>0$ such that for any choice of $\eps\in (0,\es(T))$, 
  the solution of (\ref{0eps}) satisfies
  \be{5.1}
	\int_{\{\ueps(\cdot,t)\ge 1\}} \ueps(\cdot,t)\ln \ueps(\cdot,t) \le C(T)
	\qquad \mbox{for all } t\in (0,\hatt)
  \ee
  and
  \be{5.2}
	\io \deps \frac{w_{\eps x}^2(\cdot,t)}{g(\weps(\cdot,t))} \le C(T)
	\qquad \mbox{for all } t\in (0,\hatt),
  \ee
  and such that moreover
  \be{5.3}
	\int_0^{\hatt} \io \deps \frac{u_{\eps x}^2}{\ueps} \le C(T)
  \ee
  and
  \be{5.4}
	\int_0^{\hatt} \int_{\ueps(\cdot,t)\ge 1\}} \ueps \ln\ueps \cdot f_-(\cdot,\ueps,\weps) \le C(T)
  \ee
  as well as
  \be{5.5}
	\int_0^{\hatt} \io \deps \frac{\ueps}{1+\eeps\ueps} \frac{g'(\weps)}{g(\weps)} w_{\eps x}^2 \le C(T),
  \ee
  where with $\tme$ as in Lemma \ref{lem_loc} we have set $\hatt:=\min\{T,\tme\}$.
\end{lem}
\proof
  We add the inequalities provided by Lemma \ref{lem3} and Lemma \ref{lem4} to see on dropping a nonnegative summand on the 
  right that for all $\eps\in (0,\eps_0)$ we have
  \bea{5.7}
	& & \hspace*{-20mm}
	\frac{d}{dt} \bigg\{ \io \ueps \ln \ueps + \frac{1}{2} \io \deps \frac{w_{\eps x}^2}{g(\weps)} \bigg\}
	+ \frac{1}{2} \io \deps \frac{u_{\eps x}^2}{u_\eps}
	+ \frac{1}{2} \io \deps \frac{\ueps}{1+\eeps\ueps} \frac{g'(\weps)}{g(\weps)} w_{\eps x}^2 \nn\\
	&\le& c_1
	+ \io \ueps \ln \ueps \cdot f(x,\ueps,\weps)
	+ \frac{\eps K_1}{2} \io \frac{w_{\eps x}^2}{\sqrt{g(\weps)}^3}
	\qquad \mbox{for all } t\in (0,\tme)
  \eea 
  with $c_1\equiv c_1(T):=\Big(\rho(M)+\frac{K_1}{2}\Big) \cdot \Big\{ \io u_0 \Big\} \cdot e^{\rho(M)T}$.
  Here we split $f=f_+ - f_-$ and
  \bea{5.8}
	\hspace*{-5mm}
	\io \ueps \ln \ueps \cdot f(x,\ueps,\weps)
	&=& \int_{\{\ueps<1\}} \ueps\ln\ueps \cdot f_+(x,\ueps,\weps)
	- \int_{\{\ueps<1\}} \ueps\ln\ueps \cdot f_-(x,\ueps,\weps) \nn\\
	& & + \int_{\{\ueps\ge 1\}} \ueps\ln\ueps \cdot f_+(x,\ueps,\weps)
	- \int_{\{\ueps\ge 1\}} \ueps\ln\ueps \cdot f_-(x,\ueps,\weps)
  \eea
  for $t\in (0,\tme)$, where clearly
  \be{5.9}
	\int_{\{\ueps<1\}} \ueps\ln\ueps \cdot f_+(x,\ueps,\weps) \le 0
	\qquad \mbox{for all } t\in (0,\tme),
  \ee
  and where using that
  \be{5.99}
	\xi\ln\xi \ge -\frac{1}{e}
	\qquad \mbox{for all } \xi>0,
  \ee
  we see that
  \be{5.10}
	- \int_{\{\ueps<1\}} \ueps\ln\ueps \cdot f_-(x,\ueps,\weps)
	\le \frac{|\Omega|}{e} c_2
	\qquad \mbox{for all } t\in (0,\tme)
  \ee
  with
  \bas
	c_2:=\max_{(x,u,w)\in \bar\Omega\times [0,1]\times [0,M]} f_-(x,u,w)
  \eas
  being finite by continuity of $f$.
  Since
  \be{5.999}
	f_+(\cdot,\ueps,\weps) \le \rho(M)
	\qquad \mbox{in } \Omega\times (0,\tme)
  \ee
  by (\ref{f_upper}) and Lemma \ref{lem1}, again relying on (\ref{5.99}) we see that
  \bas
	\int_{\{\ueps\ge 1\}} \ueps\ln\ueps\cdot f_+(x,\ueps,\weps)
	&\le& \rho(M) \int_{\{\ueps\ge 1\}} \ueps\ln \ueps \\
	&=& \rho(M) \cdot \bigg\{ \io \ueps\ln\ueps - \int_{\{\ueps<1\}} \ueps\ln\ueps \bigg\} \\
	&\le& \rho(M) \io \ueps\ln\ueps 
	+ \frac{\rho(M) |\Omega|}{e}
	\qquad \mbox{for all } t\in (0,\tme),
  \eas
  so that from (\ref{5.8}), (\ref{5.9}) and (\ref{5.10}) we infer that
  \be{5.11}
	\io \ueps\ln \ueps \cdot f(x,\ueps,\weps)
	\le \rho(M) \io \ueps\ln\ueps
	- \int_{\{\ueps\ge 1\}} \ueps\ln\ueps \cdot f_-(x,\ueps,\weps) + c_3
	\qquad \mbox{for all } t\in (0,\tme)
   \ee
  with $c_3:=\frac{|\Omega|}{e} \cdot (c_2+\rho(M))$.\\
  Next, the rightmost summand in (\ref{5.7}) can be estimated using Lemma \ref{lem1} along with the defining properties
  of $(\deps)_{\eps\in (0,\eps_0)}$ and $(\eeps)_{\eps\in (0,\eps_0)}$:
  Indeed, given $T>0$ we may use (\ref{eta_div}) to fix $\es=\es(T)\in (0,\eps_0)$ small enough such that
  with $\Gamma$ as in (\ref{g_upper}) we have
  \bas
	T \le \frac{\eeps}{\Gamma} \cdot \ln \frac{1}{\sqrt{\deleps}}
	\qquad \mbox{for all } \eps\in (0,\es),
  \eas
  which implies that for any such $\eps$,
  \bas
	e^{-\frac{\Gamma}{\eeps}t} 
	\ge e^{-\ln \frac{1}{\sqrt{\deleps}}}
	= \sqrt{\deleps}
	\qquad \mbox{for all } t\in (0,T)
  \eas
  and hence, by Lemma \ref{lem1},
  \bas
	M \ge \weps(x,t) 
	\ge \sqrt{\deleps} \cdot e^{-\frac{\Gamma}{\eeps}t} 
	\ge \deleps
	\qquad \mbox{for all $x\in\Omega$ and } t\in (0,\hatt).
  \eas
  Therefore, (\ref{g_pos}) applies so as to ensure that
  \bas
	g(\weps) \ge \eps
	\qquad \mbox{for all $x\in\Omega$ and } t\in (0,\hatt),
  \eas
  so that the term in question satisfies
  \bea{5.12}
	\frac{\eps K_1}{2} \io \frac{w_{\eps x}^2}{\sqrt{g(\weps)}^3}
	&\le& \frac{\sqrt{\eps} K_1}{2} \io \frac{w_{\eps x}^2}{g(\weps)} \nn\\
	&\le& \frac{K_1}{2} \io \deps \frac{w_{\eps x}^2}{g(\weps)}
	\qquad \mbox{for all } t\in (0,\hatt),
  \eea
  because $\deps\ge\sqrt{\eps}$ in $\Omega$ by (\ref{d_pos}).
  Together with (\ref{5.11}) and (\ref{5.7}), this shows that
  \bas
	\yeps(t):=\io \ueps(\cdot,t)\ln \ueps(\cdot,t)
	+ \frac{1}{2} \io \deps \frac{w_{\eps x}^2(\cdot,t)}{g(\weps(\cdot,t))},
	\qquad t\in [0,\tme),
  \eas
  and
  \bas
	\heps(t)
	&:=&
	\frac{1}{2} \io \deps \frac{u_{\eps x}^2(\cdot,t)}{\ueps(\cdot,t)}
	+ \frac{1}{2} \io \deps \frac{\ueps(\cdot,t)}{1+\eeps\ueps(\cdot,t)} \frac{g'(\weps(\cdot,t))}{g(\weps(\cdot,t))}
		w_{\eps x}^2(\cdot,t) \\
	& & + \int_{\{\ueps(\cdot,t)\ge 1\}} \ueps(\cdot,t)\ln \ueps(\cdot,t) \cdot f_-(\cdot,\ueps(\cdot,t),\weps(\cdot,t)),
	\qquad t\in (0,\tme),
  \eas
  have the property that
  \bea{5.13}
	\yeps'(t) + \heps(t)
	&\le& c_1+c_3 + \rho(M) \io \ueps\ln\ueps + \frac{K_1}{2} \io \deps \frac{w_{\eps x}^2}{g(\weps)} \nn\\
	&=& c_1+c_3 + \rho(M) \cdot \bigg\{ \yeps(t) - \frac{1}{2} \io \deps \frac{w_{\eps x}^2}{g(\weps)} \bigg\}
	+ K_1 \cdot \bigg\{ \yeps(t) - \io \ueps \ln \ueps \bigg\} \nn\\
	&\le& c_4 + c_5 \yeps(t)
	\qquad \mbox{for all } t\in (0,\hatt)
  \eea
  with $c_4:=c_1+c_3+\frac{K_1|\Omega|}{e}$ and $c_5:=\rho(M)+K_1$, where we again have used (\ref{5.99}).\\
  Now by nonnegativity of $\heps$ and (\ref{d_pos}), an integration of (\ref{5.13}) firstly yields
  \bea{5.14}
	\yeps(t)
	&\le& \yeps(0) e^{c_5 t}
	+ \frac{c_4}{c_5} e^{c_5 t} \Big(1-e^{-c_5 t}\Big) \nn\\
	&\le& c_6:=\bigg\{ \io u_0\ln u_0 + \frac{1}{2} \Big(\|d\|_{L^\infty(\Omega)}+1\Big) \cdot
	\sup_{\eps\in (0,\eps_0)} \io \frac{w_{0x}^2}{g(w_0+\deleps)} + \frac{c_4}{c_5} \bigg\} \cdot e^{c_5 T}
  \eea
  for all $t\in [0,\hatt)$ and $\eps\in (0,\es)$,
  where we note that $c_6$ is finite according to (\ref{init}), because $\deleps\to 0$ as $\eps\searrow 0$, and because due to
  the fact that $g'\ge 0$ on $[0,M]$, as asserted by (\ref{g}), Beppo Levi's theorem warrants that as $k\to\infty$ we have
  \bas
	\io \frac{w_{0 x}^2}{g(w_0+\frac{1}{k})} \nearrow \io \frac{w_{0x}^2}{g(w_0)} < \infty.
  \eas
  Once more in view of (\ref{5.99}), this entails both (\ref{5.1}) and (\ref{5.2}) with some suitably large $C(T)>0$,
  whereas another integration of (\ref{5.13}), this time making use of (\ref{5.14}), shows that
  \bea{5.15}
	\int_0^{\hatt} \heps(t)dt
	&\le& \yeps(0) - \yeps(\hatt) + c_4 \hatt + c_5 \int_0^{\hatt} \yeps(t)dt \nn\\
	&\le& c_7:=c_6 + \frac{|\Omega|}{e} + c_4 T + c_5 c_6 T,
  \eea
  and that hence (\ref{5.2})-(\ref{5.5}) hold with some possibly enlarged $C(T)$.
\qed
Due to the boundedness property (\ref{1.1}) of $\weps$, (\ref{g}) and (\ref{gamma}), from (\ref{5.2}) and
(\ref{5.5}) we particularly obtain corresponding estimates for integrals no longer containing
$\frac{1}{g(\weps)}$ and $\frac{g'(\weps)}{g(\weps)}$.
\begin{cor}\label{cor55}
  Suppose that $T>0$, and let $\es(T)\in (0,\eps_0)$ be as given by Lemma \ref{lem5}.
  Then there exists $C(T)>0$ with the property that for all $\eps\in (0,\es(T))$,
  \be{55.1}
	\io \deps w_{\eps x}^2(\cdot,t) \le C(T)
	\qquad \mbox{for all } t\in (0,\hatt)
  \ee
  and
  \be{55.2}
	\int_0^{\hatt} \io \deps \frac{\ueps}{1+\eeps\ueps} w_{\eps x}^2 \le C(T),
  \ee
  where again $\hatt:=\min\{T,\tme\}$.
\end{cor}
\proof
  Since $g(\weps) \le g(M)$ in $\Omega\times (0,\tme)$ by Lemma \ref{lem1} and (\ref{g}), we immediately obtain 
  (\ref{55.1}) from (\ref{5.2}).
  As furthermore (\ref{gamma}) warrants that
  \bas
	\frac{g'(\weps)}{g(\weps)} \ge \gamma>0
	\qquad \mbox{in } \Omega\times (0,\tme)
  \eas
  by Lemma \ref{lem1}, we also infer from (\ref{5.5}) that (\ref{55.2}) is valid with some adequately large
  $C(T)>0$.
\qed
As one consequence of (\ref{55.1}) when combined with the boundedness of $\frac{\ueps}{(1+\eeps\ueps)^2}$ and the uniform
positivity of $\deps$, we can infer that in fact our approximate solutions cannot blow up in finite time:
\begin{lem}\label{lem_global}
  For each $\eps\in (0,\eps_0)$, the solution of (\ref{0eps}) is global in time; that is, in Lemma \ref{lem_loc}
  we have $\tme=\infty$.
\end{lem}
\proof
  Assuming on the contrary that $\tme$ be finite, combining (\ref{1.2}) with (\ref{g}) we first obtain that then
  there would exist $c_1>0$ such that
  \be{gl1}
	\frac{1}{g(\weps)} \le c_1
	\qquad \mbox{in } \Omega\times (0,\tme).
  \ee
  Moreover, as $\deps>0$ in $\bar\Omega$ by (\ref{d_pos}), Corollary \ref{cor55} and Lemma \ref{lem1}
  would yield $c_2>0$ fulfilling
  \be{gl2}
	\|\weps(\cdot,t)\|_{W^{1,2}(\Omega)} \le c_2
	\qquad \mbox{for all } t\in (0,\tme).
  \ee
  In particular, the latter along with (\ref{d_pos}) and the fact that 
  $\frac{\xi}{(1+\eeps\xi)^2} \le \frac{1}{4\eeps}$ for all
  $\xi\ge 0$ ensures that the cross-diffusive flux in the first equation in (\ref{0eps}) satisfies
  \bas
	\Big\|\deps \frac{\ueps(\cdot,t)}{(1+\eeps\ueps(\cdot,t))^2} w_{\eps x}(\cdot,t) \Big\|_{L^2(\Omega)}
	\le \Big(\|d\|_{L^\infty(\Omega)}+1 \Big) \cdot \frac{1}{4\eeps} \cdot c_2
	\qquad \mbox{for all } t\in (0,\tme).
  \eas
  Since furthermore, by (\ref{f_upper}) and again Lemma \ref{lem1}, 
  \bas
	f(\cdot,\ueps,\weps) \le \rho(M)
	\qquad \mbox{in } \Omega\times (0,\tme),
  \eas
  a standard argument based on the smoothing properties of the non-degenerate linear semigroup
  $(e^{t(\deps \cdot)_{xx}})_{t\ge 0}$ (cf.~e.g.~the reasoning in \cite[Lemma 3.2]{BBTW}) provides $c_3>0$ such that
  \bas
	\|\ueps(\cdot,t)\|_{L^\infty(\Omega)} \le c_3
	\qquad \mbox{for all } t\in (0,\tme).
  \eas
  In view of the extensibility criterion (\ref{ext}), together with (\ref{gl1}) and (\ref{gl2}) this shows that
  our assumption $\tme<\infty$ was absurd.
\qed
\mysection{Weak precompactness properties of $\ueps f(x,\ueps,\weps)$ and $\ueps g(\weps)$ in $L^1$}\label{precomp}
In appropriately passing to the limit in the zero-order integrals appearing in the
respective weak formulations of (\ref{0eps}), we shall make essential use of 
two compactness properties of the solutions thereof which appear to go beyond trivial implications of the bounds
provided by Lemma \ref{lem5}.
As a preparation for our arguments in this respect, let us state the following observation on a lower bound for
all possible values of $u\ge 0$ at which $u\cdot f_-(x,u,w)$ may become large for some $x\in\bar\Omega$ and
$w\in [0,M]$. This will be used in Lemma \ref{lem6} to assert that for arbitrarily large $\kappa>0$ one can pick
$N>0$ in such a way that whenever $\ueps f_-(\cdot\ueps,\weps) \ge N$, we know that $\ueps\ge \kappa$.
\begin{lem}\label{lem7}
  With $M>0$ as in (\ref{M}), let
  \be{7.1}
	\set(N):= \Big\{ u\ge 0 \ \Big| \ u\cdot f_-(x,u,w) \ge N
	\mbox{ for some $x\in\bar\Omega$ and $w\in [0,M]$} \Big\}
  \ee
  and
  \be{7.2}
	\kappa(N):=\left\{ \begin{array}{ll}
	\inf \set(N)
	\qquad & \mbox{if } \set(N)\ne\emptyset, \\[1mm]
	+\infty & \mbox{else}
	\end{array} \right.
  \ee
  for $N\in\N$.
  Then 
  \be{7.3}
	\limsup_{N\to\infty} \kappa(N)=+\infty.
  \ee
\end{lem}
\proof
  If (\ref{7.3}) was false, then there would exist $N_0\in\N$ such that for all $N\ge N_0$ we would have 
  $\set(N)\ne\emptyset$ and $\kappa(N)<c_1$ with some $c_1>0$.
  By definition of $\set(N)$ and $\kappa(N)$, this would mean that we could find $(x_N)_{N\ge N_0} \subset\bar\Omega$,
  $(u_N)_{N\ge N_0} \subset [0,\infty)$ and $(w_N)_{N\ge N_0} \subset [0,M]$ fulfilling
  \be{7.4}
	u_N \cdot f_-(x_N,u_N,w_N) \ge N
	\qquad \mbox{for all } N\ge N_0
  \ee
  and
  \bas
	u_N \le c_1
	\qquad \mbox{for all } N\ge N_0,
  \eas
  where passing to a subsequence if necessary we may assume that as $N\to\infty$ we have $x_N\to x_\infty$,
  $u_N\to u_\infty$ and $w_N\to w_\infty$ with some $x_\infty\in\bar\Omega$, $u_\infty \in [0,c_1]$ and
  $w_\infty \in [0,M]$.
  By continuity of $f_-$, however, this would imply that 
  \bas
	u_N \cdot f_-(x_N,u_N,w_N) \to u_\infty \cdot f_-(x_\infty,u_\infty,w_\infty)
	\qquad \mbox{as } N\to\infty
  \eas
  and thereby contradict (\ref{7.4}).
\qed
Making use of the latter, by means of the Dunford-Pettis theorem we can now establish suitable compactness
properties of the rightmost summands in the first two equations in (\ref{0eps}).
\begin{lem}\label{lem6}
  Let $T>0$. Then with $\es(T)\in (0,\eps_0)$ as in Lemma \ref{lem5},
  \be{6.1}
	\Big(\ueps f(\cdot,\ueps,\weps)\Big)_{\eps\in (0,\es(T))}
	\quad \mbox{is relatively compact with respect to the weak topology in } L^1(\Omega\times (0,T)),
  \ee
  and moreover
  \be{6.01}
	\Big(\ueps g(\weps)\Big)_{\eps\in (0,\es(T))}
	\quad \mbox{is relatively compact with respect to the weak topology in } L^1(\Omega\times (0,T)).
  \ee
\end{lem}
\proof
  According to Lemma \ref{lem5}, we can fix positive constants $c_1$ and $c_2$ such that
  \be{6.2}
	\int_{\{\ueps(\cdot,t)\ge 1\}} \ueps(\cdot,t)\ln \ueps(\cdot,t) \le c_1
	\qquad \mbox{for all } t\in (0,T)
  \ee
  and
  \be{6.3}
	\int_0^T \int_{\ueps(\cdot,t)\ge 1\}} \ueps \ln \ueps \cdot f_- (\cdot,\ueps,\weps) \le c_2
  \ee
  whenever $\eps\in (0,\es(T))$.
  Aiming at an application of the Dunford-Pettis theorem, given $\mu>0$ we first fix an integer $N\ge 1$ suitably large
  such that
  \be{6.4}
	\frac{c_1 \rho(M) T}{\ln N} < \frac{\mu}{4}
  \ee
  and
  \be{6.44}
	\frac{c_1 g(M)T}{\ln N} < \frac{\mu}{2},
  \ee
  and such that with $\kappa(N)$ as defined in Lemma \ref{lem7} we have $\kappa(N)>1$ and
  \be{6.5}
	\frac{c_2}{\ln \kappa(N)} <\frac{\mu}{4},
  \ee
  where the latter is possible due to the outcome of Lemma \ref{lem7}.
  Thereafter, we choose $\iota>0$ small enough fulfilling
  \be{6.6}
	\rho(M) N \iota < \frac{\mu}{4}
  \ee
  and
  \be{6.7}
	N\iota <\frac{\mu}{4}
  \ee
  as well as
  \be{6.76}
	g(M) N\iota < \frac{\mu}{2},
  \ee
  and fix an arbitrary measurable set $E\subset \Omega\times (0,T)$ satisfying $|E|<\iota$. Then decomposing
  \be{6.77}
	\int\int_E \Big|\ueps f(\cdot,\ueps,\weps)\Big|
	= \int\int_E \ueps f_+(\cdot,\ueps,\weps)
	+ \int\int_I \ueps f_-(\cdot,\ueps,\weps),
  \ee
  by combining (\ref{6.2}) with Lemma \ref{lem1}, (\ref{f_upper}) and (\ref{6.6}) we can estimate 
  \bea{6.8}
	\int_E \ueps f_+(\cdot\ueps,\weps)
	&=& \int\int_{E\cap \{\ueps<N\}} \ueps f_+(\cdot,\ueps,\weps)
	+ \int_{E\cap \{\ueps\ge N\}} \ueps f_+(\cdot,\ueps,\weps) \nn\\
	&\le& \rho(M) \in\int_{E\cap \{\ueps<N\}} \ueps
	+ \frac{\rho(M)}{\ln N} \int\int_{E\cap \{\ueps\ge N\}} \ueps\ln \ueps \nn\\
	&\le& \rho(M) \cdot N|E|
	+ \frac{\rho(M)}{\ln N} \int_0^T \int_{\{\ueps(\cdot,t)\ge 1\}} \ueps \ln \ueps \nn\\
	&\le& \rho(M) N\iota
	+ \frac{\rho(M)}{\ln N} \cdot c_1 T \nn\\
	&<& \frac{\mu}{4}+\frac{\mu}{4}=\frac{\mu}{2}
	\qquad \mbox{for all } \eps\in (0,\es(T)).
  \eea
  Likewise, relying on (\ref{6.7}) we see that
  \bea{6.9}
	\int\int_E \ueps f_-(\cdot,\ueps,\weps)
	&=& \int\int_{E\cap \{\ueps f_-(\cdot,\ueps,\weps)<N\}} \ueps f_-(\cdot,\ueps,\weps)
	+ \int\int_{E\cap \{\ueps f_-(\cdot,\ueps,\weps)\ge N\}} \ueps f_-(\cdot,\ueps,\weps) \nn\\
	&\le& N|E|
	+ \int\int_{E\cap \{\ueps f_-(\cdot,\ueps,\weps)\ge N\}} \ueps f_-(\cdot,\ueps,\weps) \nn\\
	&<& \frac{\mu}{4} 
	+ \int\int_{E\cap \{\ueps f_-(\cdot,\ueps,\weps)\ge N\}} \ueps f_-(\cdot,\ueps,\weps)
	\qquad \mbox{for all } \eps\in (0,\es(T)),
  \eea
  and in order to appropriately control the last summand herein we recall the definition (\ref{7.2}) of $\kappa(N)$ 
  to observe that whenever $\ueps(x,t) f_-(x,\ueps(x,t),\weps(x,t))\ge N$ for some $\eps\in (0,\eps_0)$, $x\in\bar\Omega$
  and $t\ge 0$, we necessarily must have $\ueps(x,t) \ge \kappa(N)$.
  Consequently, $E\cap \{\ueps f_-(\cdot\ueps,\weps) \ge N\} \subset E\cap \{\ueps\ge \kappa(N)\}$,
  so that (\ref{6.3}) and (\ref{6.5}) become applicable so as to guarantee that
  \bas
	\int\int_{E\cap \{\ueps f_-(\cdot,\ueps,\weps)\ge N\}} \ueps f_-(\cdot,\ueps,\weps)
	&\le& \int\int_{E\cap \{\ueps \ge \kappa(N)}\} \ueps f_-(\cdot,\ueps,\weps) \\
	&\le& \frac{1}{\ln \kappa(N)} \int\int_{E\cap \{\ueps\ge \kappa(N)\}} \ueps\ln\ueps \cdot f_-(\cdot,\ueps,\weps) \\
	&\le& \frac{1}{\ln\kappa(N)} \cdot c_2 \\
	&<& \frac{\mu}{4}
	\qquad \mbox{for all } \eps\in (0,\es(T)),
  \eas
  which along with (\ref{6.77}), (\ref{6.8}) and (\ref{6.9}) shows that for any such $E$ we have
  \be{6.10}
	\int\int_E \Big|\ueps f	(\cdot,\ueps,\weps)\Big| < \mu
	\qquad \mbox{for all } \eps\in (0,\es(T)).
  \ee
  Similarly, using Lemma \ref{lem1} together with (\ref{g}) we obtain
  \bas
	\int\int_E \Big|\ueps g(\weps)\Big|
	&=& \int\int_{E\cap \{\ueps<N\}} \ueps g(\weps)
	+ \int\int_{E\cap \{\ueps \ge N\}} \ueps g(\weps) \\
	&\le& g(M) \int\int_{E\cap\{\ueps<N\}} \ueps
	+ g(M) \int\int_{E\cap \{\ueps\ge N\}} \ueps \\
	&\le& g(M) \cdot N|E|
	+ \frac{g(M)}{\ln N} \int_0^T \int_{\{\ueps(\cdot,t)\ge 1\}} \ueps \\
	&\le& g(M) N\iota
	+ \frac{g(M)}{\ln N} \cdot c_1 T \\
	&<& \frac{\mu}{2}+\frac{\mu}{2}=\mu
	\qquad \mbox{for all } \eps\in (0,\es(T)),
  \eas
  because of (\ref{6.44}) and (\ref{6.76}).
  By means of the Dunford-Pettis compactness criterion, from this we infer that (\ref{6.01}) holds, and that
  (\ref{6.1}) is a consequence of (\ref{6.10}).
\qed
\mysection{Regularity properties of $\sqrt{\deps} \ueps$}\label{reg-prop}
In order to further prepare our limit procedure, especially with regard to pointwise convergence of $\ueps$ and
of convergence in
the cross-diffusive flux term $\deps \frac{\ueps}{(1+\eeps\ueps)^2} w_{\eps x}$ in (\ref{0eps}), we next plan to combine
the weak compactness feature of the part $\sqrt{\deps}w_{\eps x}$ thereof,
as naturally implied by Corollary \ref{cor55}, 
by an appropriate result on convergence in the complementary factor 
$\sqrt{\deps}\frac{\ueps}{(1+\eeps\ueps)^2}$ in a strong $L^2$ topology.
To achieve this in Lemma \ref{lem12} by using underway an argument based on the Aubin-Lions lemma,
let us suitably interpolate between the inequalities in (\ref{2.2}) and (\ref{5.3}) to derive the following
spatio-temporal estimates for the quantity $\sqrt{\deps}\ueps$ forming the core of the factor in question.
\begin{lem}\label{lem11}
  Let $T>0$ and $\es(T)\in (0,\eps_0)$ be as in Lemma \ref{lem5}.
  Then there exists $C(T)>0$ such that for all $\eps\in (0,\es(T))$,
  \be{11.1}
	\int_0^T \Big\| (\sqrt{\deps} \ueps(\cdot,t))_x\Big\|_{L^1(\Omega)}^2 dt \le C(T)
  \ee
  and
  \be{11.01}
	\int_0^T \Big\|\sqrt{\deps} \ueps(\cdot,t)\Big\|_{L^\infty(\Omega)}^2 dt \le C(T)
  \ee
  as well as
  \be{11.2}
	\int_0^T \io \sqrt{\deps}^3 \ueps^3 \le C(T).
  \ee
\end{lem}
\proof
  According to Lemma \ref{lem2} and Lemma \ref{lem5}, there exist $c_1=c_1(T)>0$ and $c_2=c_2(T)>0$ such that
  for any $\eps\in (0,\eps_0)$ we have
  \be{11.3}
	\io \ueps \le c_1
	\qquad \mbox{for all } t\in (0,T),
  \ee
  and that
  \be{11.4}
	\int_0^T \io \deps \frac{u_{\eps x}^2}{\ueps} \le c_2
  \ee
  whenever $\eps\in (0,\es(T))$.
  Since
  \bas
	\Big|(\sqrt{\deps}\ueps)_x\Big|
	&=& \Big|\sqrt{\deps} u_{\eps x} + \frac{d_{\eps x}}{2\sqrt{\deps}} \ueps \Big| \\
	&\le& \sqrt{\deps} |u_{\eps x}|
	+ \frac{\sqrt{K_1}}{2} \ueps
	\qquad \mbox{in } \Omega\times (0,\infty)
  \eas
  due to (\ref{K1}), by the Cauchy-Schwarz inequality these estimates imply that
  \bas
	\int_0^T \Big\| (\sqrt{\deps}\ueps(\cdot,t))_x\Big\|_{L^1(\Omega)}^2 dt
	&\le& \int_0^T \bigg\{ \io \sqrt{\deps} |u_{\eps x}| \bigg\}^2
	+ \frac{\sqrt{K_1}}{2} \int_0^T \bigg\{ \io \ueps\bigg\}^2 \\
	&\le& \int_0^T \bigg\{ \io \deps \frac{u_{\eps x}^2}{\ueps} \bigg\} \cdot \bigg\{ \io \ueps \bigg\}
	+ \frac{\sqrt{K_1}}{2} \int_0^T \bigg\{ \io \ueps\bigg\}^2 \\
	&\le& c_3\equiv c_3(T):=c_1 c_2 + \frac{\sqrt{K_1}}{2} c_1^2 T
	\qquad \mbox{for all } \eps\in (0,\es(T))
  \eas
  and thereby proves (\ref{11.1}), whereupon (\ref{11.01}) follows from Lemma \ref{lem2} and
  the fact that $W^{1,1}(\Omega) \hra L^\infty(\Omega)$.
  As the Gagliardo-Nirenberg inequality provides $c_4>0$ such that
  \bas
	\|\varphi\|_{L^3(\Omega)}^3 \le c_4 \|\varphi_x\|_{L^1(\Omega)}^2 \|\varphi\|_{L^1(\Omega)}
	+ c_4\|\varphi\|_{L^1(\Omega)}^3
	\qquad \mbox{for all } \varphi \in W^{1,1}(\Omega),
  \eas
  in view of (\ref{d_pos}) this furthermore entails that with $c_5:=\sqrt{\|d\|_{L^\infty(\Omega)}+1}$ we have
  \bas
	\int_0^T \io \sqrt{\deps}^3 \ueps^3
	&=& \int_0^T \Big\|\sqrt{\deps}\ueps(\cdot,t)\Big\|_{L^3(\Omega)}^3 dt \\
	&\le& c_4 \int_0^T \Big\|(\sqrt{\deps}\ueps(\cdot,t))_x\Big\|_{L^1(\Omega)}^2
	\cdot \Big\| \sqrt{\deps}\ueps(\cdot,t)\Big\|_{L^1(\Omega)} dt 
	+ c_4 \int_0^T \Big\|\sqrt{\deps}\ueps(\cdot,t)\Big\|_{L^1(\Omega)}^3 dt \\
	&\le& c_4 c_5 \big\|\ueps\big\|_{L^\infty((0,T);L^1(\Omega))}
	\int_0^T \Big\|(\sqrt{\deps}\ueps(\cdot,t))_x\Big\|_{L^1(\Omega)}^2 dt
	+ c_4 c_5^3 \int_0^T \Big\|\ueps(\cdot,t)\Big\|_{L^1(\Omega)}^3 dt \\[2mm]
	&\le& c_4 c_5 c_1 c_3 + c_4 c_5^3 c_1^3 T
	\qquad \mbox{for all } \eps\in (0,\es(T)),
  \eas
  and that thus also (\ref{11.2}) holds.
\qed
\mysection{Regularity in time}\label{reg-time}
As a final preparation for our first subsequence extraction procedure, we combine our previously gained estimates
to obtain some regularity properties involving time derivatives of the solution components $\ueps$ and $\weps$.
\begin{lem}\label{lem31}
  Let $T>0$ and $\es(T)\in (0,\eps_0)$ be as in Lemma \ref{lem5}.
  Then there exists $C(T)>0$ such that
  \be{31.1}
	\int_0^T \Big\|\partial_t \sqrt{\deps(\ueps(\cdot,t)+1)} \Big\|_{(W^{1,3}(\Omega))^\star} dt \le C(T)
	\qquad \mbox{for all } \eps\in (0,\es(T)).
  \ee
\end{lem}
\proof
  For fixed $t>0$ and $\psi\in C^1(\bar\Omega)$, from the first equation in (\ref{0eps}) we obtain that
  \bea{31.2}
	\io \partial_t \sqrt{\deps(\ueps(\cdot,t)+1)}{\psi} 
	&=& - \frac{1}{2} \io \Big( \sqrt{\deps} \frac{1}{\sqrt{\ueps+1}} \psi\Big)_x \cdot (\deps\ueps)_x \nn\\
	& & + \frac{1}{2} \io \Big( \sqrt{\deps} \frac{1}{\sqrt{\ueps+1}} \psi\Big)_x \cdot 
		\deps \frac{\ueps}{(1+\eeps\ueps)^2} w_{\eps x} \nn\\
	& & + \frac{1}{2} \io \sqrt{\deps} \frac{1}{\sqrt{\ueps+1}} \ueps f(\cdot,\ueps,\weps) \psi \nn\\[2mm]
	&=& \frac{1}{4} \io \sqrt{\deps}^3 \frac{1}{\sqrt{\ueps+1}^3} u_{\eps x}^2 \psi
	+ \frac{1}{4} \io \sqrt{\deps} d_{\eps x} \frac{\ueps}{\sqrt{\ueps+1}^3} u_{\eps x} \psi \nn\\
	& & - \frac{1}{4} \io \sqrt{\deps} d_{\eps x} \frac{1}{\sqrt{\ueps+1}} u_{\eps x} \psi
	- \frac{1}{4} \io \frac{d_{\eps x}^2}{\sqrt{\deps}} \frac{\ueps}{\sqrt{\ueps+1}} \psi \nn\\
	& & - \frac{1}{2} \io \sqrt{\deps}^3 \frac{1}{\sqrt{\ueps+1}} u_{\eps x} \psi_x
	- \frac{1}{2} \io \sqrt{\deps} d_{\eps x} \frac{\ueps}{\sqrt{\ueps+1}} \psi_x \nn\\
	& & - \frac{1}{4} \io \sqrt{\deps}^3 \frac{\ueps}{\sqrt{\ueps+1}^3 (1+\eeps\ueps)^2} u_{\eps x} w_{\eps x} \psi \nn\\
	& & + \frac{1}{4} \io \sqrt{\deps} d_{\eps x} \frac{\ueps}{\sqrt{\ueps+1}(1+\eeps\ueps)^2} w_{\eps x} \psi \nn\\
	& & + \frac{1}{2} \io \sqrt{\deps}^3 \frac{\ueps}{\sqrt{\ueps+1} (1+\eeps\ueps)^2} w_{\eps x} \psi_x \nn\\
	& & + \frac{1}{2} \io \sqrt{\deps} \frac{\ueps}{\sqrt{\ueps+1}} f(\cdot,\ueps,\weps) \psi
	\qquad \mbox{for all } \eps\in (0,\eps_0).
  \eea
  Here using (\ref{d_pos}) and Young's inequality, we see that
  \bea{31.3}
	\bigg| \frac{1}{4} \io \sqrt{\deps}^3 \frac{1}{\sqrt{\ueps+1}^3} u_{\eps x}^2 \psi \bigg|
	&\le& \frac{\sqrt{c_1}}{4} \bigg\{ \io \deps \frac{u_{\eps x}^2}{\ueps} \bigg\}^\frac{1}{2} 
		\|\psi\|_{L^\infty(\Omega)} \nn\\
	&\le& \frac{\sqrt{c_1}}{8} \bigg\{ \io \deps \frac{u_{\eps x}^2}{\ueps} + 1 \bigg\} 
		\|\psi\|_{L^\infty(\Omega)}
  \eea
  with $c_1:=\|d\|_{L^\infty(\Omega)}+1$, and similarly, 
  \bea{31.4}
	\bigg| \frac{1}{4} \io \sqrt{\deps} d_{\eps x} \frac{\ueps}{\sqrt{\ueps+1}^3} u_{\eps x} \psi \bigg|
	&\le& \frac{1}{4} \bigg\{ \io \deps \frac{u_{\eps x}^2}{\ueps} \bigg\}^\frac{1}{2}
	\cdot \bigg\{ \io d_{\eps x}^2 \frac{\ueps^3}{(\ueps+1)^3} \psi^2 \bigg\}^\frac{1}{2} \nn\\
	&\le& \frac{\sqrt{K_1 c_1}}{4} \bigg\{ \io \deps \frac{u_{\eps x}^2}{\ueps} \bigg\}^\frac{1}{2}
		\|\psi\|_{L^2(\Omega)} \nn\\
	&\le& \frac{\sqrt{K_1 c_1}}{8} \bigg\{ \io \deps \frac{u_{\eps x}^2}{\ueps} +1 \bigg\}
		\|\psi\|_{L^2(\Omega)},
  \eea
  because
  \be{31.5}
	d_{\eps x}^2 \le K_1 \deps \le K_1 c_1
	\qquad \mbox{in } \Omega
  \ee
  thanks to (\ref{K1}) and (\ref{d_pos}).
  Next, by the H\"older inequality and again due to (\ref{31.5}), (\ref{K1}), (\ref{d_pos}) and Young's inequality, 
  \bea{31.61}
	\bigg|- \frac{1}{4} \io \sqrt{\deps} d_{\eps x} \frac{1}{\sqrt{\ueps+1}} u_{\eps x} \psi \bigg|
	&\le& \frac{1}{4} \bigg\{ \io \deps \frac{u_{\eps x}^2}{\ueps} \bigg\}^\frac{1}{2}
	\cdot \bigg\{ \io d_{\eps x}^2 \frac{\ueps}{\ueps+1} \psi^2 \bigg\}^\frac{1}{2} \nn\\
	&\le& \frac{\sqrt{K_1 c_1}}{4} \bigg\{ \io \deps \frac{u_{\eps x}^2}{\ueps} \bigg\}^\frac{1}{2}
		\|\psi\|_{L^2(\Omega)} \nn\\
	&\le& \frac{\sqrt{K_1 c_1}}{8} \bigg\{ \io \deps \frac{u_{\eps x}^2}{\ueps} +1 \bigg\}
		\|\psi\|_{L^2(\Omega)} 
  \eea
  and
  \bea{31.62}
	\bigg|- \frac{1}{4} \io \frac{d_{\eps x}^2}{\sqrt{\deps}} \frac{\ueps}{\sqrt{\ueps+1}} \psi \bigg|
	&\le& \frac{1}{4} \bigg\{ \io \sqrt{\deps}^3 \ueps^3 \bigg\}^\frac{1}{6} 
	\cdot \bigg\{ \io \deps^{-\frac{9}{10}} |d_{\eps x}|^\frac{12}{5} \Big(\frac{\ueps}{\ueps+1}\Big)^\frac{3}{5}
		|\psi|^\frac{6}{5} \bigg\}^\frac{5}{6} \nn\\
	&\le& \frac{K_1}{4} \bigg\{ \io \sqrt{\deps}^3 \ueps^3 \bigg\}^\frac{1}{6}
	\cdot \bigg\{ \io \deps^\frac{3}{10} |\psi|^\frac{6}{5}\bigg\}^\frac{5}{6} \nn\\
	&\le& \frac{K_1 c_1^\frac{1}{4}}{4} \bigg\{ \io \sqrt{\deps}^3 \ueps^3 \bigg\}^\frac{1}{6}
	\|\psi\|_{L^\frac{6}{5}(\Omega)} \nn\\
	&\le& \frac{K_1 c_1^\frac{1}{4}}{24} \bigg\{ \io \sqrt{\deps}^3 \ueps^3 +1 \bigg\}
	\|\psi\|_{L^\frac{6}{5}(\Omega)} 
  \eea
  as well as
  \bea{31.71}
	\bigg|- \frac{1}{2} \io \sqrt{\deps}^3 \frac{1}{\sqrt{\ueps+1}} u_{\eps x} \psi_x \bigg|
	&\le& \frac{1}{2} \bigg\{ \io \deps \frac{u_{\eps x}^2}{\ueps} \bigg\}^\frac{1}{2}
	\cdot \bigg\{ \io \deps^2 \frac{\ueps}{\ueps+1} \psi_x^2 \bigg\}^\frac{1}{2} \nn\\
	&\le& \frac{c_1}{2} \bigg\{ \io \deps \frac{u_{\eps x}^2}{\ueps} \bigg\}^\frac{1}{2}
	\|\psi_x\|_{L^2(\Omega)} \nn\\
	&\le& \frac{c_1}{4} \bigg\{ \io \deps \frac{u_{\eps x}^2}{\ueps} +1 \bigg\}
	\|\psi_x\|_{L^2(\Omega)} 
  \eea
  and
  \bea{31.72}
	\bigg|- \frac{1}{2} \io \sqrt{\deps} d_{\eps x} \frac{\ueps}{\sqrt{\ueps+1}} \psi_x \bigg|
	&\le& \frac{1}{2} \bigg\{ \io \sqrt{\deps}^3 \ueps^3 \bigg\}^\frac{1}{6}
	\cdot \bigg\{ \deps^\frac{3}{10} |d_{\eps x}|^\frac{6}{5} \Big(\frac{\ueps}{\ueps+1}\Big)^\frac{3}{5}
		|\psi_x|^\frac{6}{5} \bigg\}^\frac{5}{6} \nn\\
	&\le& \frac{\sqrt{K_1}}{2} \bigg\{ \io \sqrt{\deps}^3 \ueps^3 \bigg\}^\frac{1}{6}
	\cdot \bigg\{ \deps^\frac{9}{10} |\psi_x|^\frac{6}{5} \bigg\}^\frac{5}{6} \nn\\
	&\le& \frac{\sqrt{K_1} c_1^\frac{3}{4}}{2} \bigg\{ \io \sqrt{\deps}^3 \ueps^3 \bigg\}^\frac{1}{6}
	\|\psi_x\|_{L^\frac{6}{5}(\Omega)} \nn\\
	&\le& \frac{\sqrt{K_1} c_1^\frac{3}{4}}{2} \bigg\{ \io \sqrt{\deps}^3 \ueps^3 +1 \bigg\}
	\|\psi_x\|_{L^\frac{6}{5}(\Omega)}.
  \eea
  Likewise, the integrals in (\ref{31.2}) stemming from the cross-diffusive interaction can be estimated according to
  \bea{31.8}
	\bigg|- \frac{1}{4} \io \sqrt{\deps}^3 \frac{\ueps}{\sqrt{\ueps+1}^3 (1+\eeps\ueps)^2} u_{\eps x} w_{\eps x} \psi
		\bigg|
	&\le& \frac{1}{4} \bigg\{ \io \deps \frac{u_{\eps x}^2}{\ueps} \bigg\}^\frac{1}{2}
	\cdot \bigg\{ \io \deps^2 \frac{\ueps^3}{(\ueps+1)^3} w_{\eps x}^2 \bigg\}^\frac{1}{2}
		\|\psi\|_{L^\infty(\Omega)} \nn\\
	&\le& \frac{\sqrt{c_1}}{4} \bigg\{ \io \deps \frac{u_{\eps x}^2}{\ueps} \bigg\}^\frac{1}{2}
	\cdot \bigg\{ \io \deps w_{\eps x}^2 \bigg\}^\frac{1}{2}
		\|\psi\|_{L^\infty(\Omega)} \nn\\
	&\le& \frac{\sqrt{c_1}}{8} \bigg\{ \io \deps \frac{u_{\eps x}^2}{\ueps} + \io \deps w_{\eps x}^2 \bigg\}
		\|\psi\|_{L^\infty(\Omega)}
  \eea
  and
  \bea{31.9}
	\bigg|\frac{1}{4} \io \sqrt{\deps} d_{\eps x} \frac{\ueps}{\sqrt{\ueps+1}(1+\eeps\ueps)^2} w_{\eps x} \psi \bigg|
	&\le& \frac{1}{4} \bigg\{ \io \deps w_{\eps x}^2 \bigg\}^\frac{1}{2}
	\cdot \bigg\{ \io d_{\eps x}^2 \frac{\ueps^2}{\ueps+1} \psi^2 \bigg\}^\frac{1}{2} \nn\\
	&\le& \frac{\sqrt{K_1 c_1}}{4} \bigg\{ \io \deps w_{\eps x}^2 \bigg\}^\frac{1}{2}
	\cdot \bigg\{ \io \ueps \bigg\}^\frac{1}{2} \|\psi\|_{L^2(\Omega)} \nn\\
	&\le& \frac{\sqrt{K_1 c_1}}{8} \bigg\{ \io \deps w_{\eps x}^2 + \io \ueps \bigg\} \|\psi\|_{L^2(\Omega)}
  \eea
  as well as
  \bea{31.10}
	\bigg|\frac{1}{2} \io \sqrt{\deps}^3 \frac{\ueps}{\sqrt{\ueps+1} (1+\eeps\ueps)^2} w_{\eps x} \psi_x \bigg| 
	&\le& \frac{1}{2} \bigg\{ \io \deps w_{\eps x}^2 \bigg\}^\frac{1}{2}
	\cdot \bigg\{ \io \deps^2 \frac{\ueps^2}{\ueps+1} \psi_x^2 \bigg\}^\frac{1}{2} \nn\\
	&\le& \frac{c_1^\frac{3}{4}}{2} \bigg\{ \io \deps w_{\eps x}^2 \bigg\}^\frac{1}{2}
	\cdot \bigg\{ \io \sqrt{\deps} \ueps \psi_x^2 \bigg\}^\frac{1}{2} \nn\\
	&\le& \frac{c_1^\frac{3}{4}}{2} \bigg\{ \io \deps w_{\eps x}^2 \bigg\}^\frac{1}{2}
	\cdot \bigg\{ \io \sqrt{\deps}^3 \ueps^3 \bigg\}^\frac{1}{6} \|\psi_x\|_{L^3(\Omega)} \nn\\
	&\le& \frac{c_1^\frac{3}{4}}{4} \bigg\{ \io \deps w_{\eps x}^2 + \io \sqrt{\deps}^3 \ueps^3 + 1 \bigg\}
	\|\psi_x\|_{L^3(\Omega)}.
  \eea
  Since finally 
  \bea{31.11}
	\bigg|\frac{1}{2} \io \sqrt{\deps} \frac{\ueps}{\sqrt{\ueps+1}} f(\cdot,\ueps,\weps) \psi \bigg| 
	\le \frac{c_1}{2} \bigg\{ \io \ueps |f(\cdot,\ueps,\weps)| \bigg\} \|\psi\|_{L^\infty(\Omega)},
  \eea
  and since in the present one-dimensional setting we have
  $W^{1,3}(\Omega) \subset W^{1,\frac{6}{5}}(\Omega) \hra L^\infty(\Omega) \subset
  L^3(\Omega)\subset L^2(\Omega)\subset L^\frac{6}{5}(\Omega)$,
  in view of the estimates implied by Lemma \ref{lem5}, Corollary \ref{cor55}, Lemma \ref{lem11}, Lemma \ref{lem2} and
  Lemma \ref{lem6} we only need to collect (\ref{31.3}), (\ref{31.4}) and (\ref{31.61})-(\ref{31.11})
  to derive (\ref{31.1}) from (\ref{31.2}).
\qed
It may be not surprising that our derivation of a corresponding property of $\weps$ is much less involved:
\begin{lem}\label{lem32}
  Let $T>0$, and let $\es(T)\in (0,\eps_0)$ be as in Lemma \ref{lem5}.
  Then one can find $C(T)>0$ such that
  \be{32.1}
	\int_0^T \Big\|\partial_t \Big(\sqrt{\deps} \weps(\cdot,t)\Big) \Big\|_{(W^{1,2}(\Omega))^\star}^3 dt \le C(T)
	\qquad \mbox{for all } \eps\in (0,\es(T)).
  \ee
\end{lem}
\proof
  For arbitrary $\psi\in C^1(\bar\Omega)$, the second equation in (\ref{0eps}) shows that
  \bea{32.2}
	\io \partial_t \Big(\sqrt{\deps} \weps(\cdot,t)\Big) \psi
	&=& - \eps \io \frac{w_{\eps x}}{\sqrt{g(\weps)}} \big(\sqrt{\deps} \psi\big)_x
	- \io \sqrt{\deps} \frac{\ueps}{1+\eta _\eps u_\eps} g(\weps)\psi \nn\\
	&=& - \frac{\eps}{2} \io \frac{d_{\eps x}}{\sqrt{\deps}} \frac{w_{\eps x}}{\sqrt{g(\weps)}} \psi
	- \eps \io \sqrt{\deps} \frac{w_{\eps x}}{\sqrt{g(\weps)}} \psi_x \nn \\
	&-& \io \sqrt{\deps} \frac{\ueps}{1+\eta _\eps u_\eps} g(\weps)\psi
  \eea
  for all $\eps\in (0,\eps_0)$, where by the Cauchy-Schwarz inequality and (\ref{K1}),
  \bea{32.3}
	\bigg|- \frac{\eps}{2} \io \frac{d_{\eps x}}{\sqrt{\deps}} \frac{w_{\eps x}}{\sqrt{g(\weps)}} \psi \bigg|
	&\le& \frac{\eps}{2} \bigg\{ \io \deps \frac{w_{\eps x}^2}{g(\weps)} \bigg\}^\frac{1}{2}
	\cdot \bigg\{ \io \frac{d_{\eps x}^2}{\deps^2} \psi^2 \bigg\}^\frac{1}{2} \nn\\
	&\le& \frac{\sqrt{K_1}\eps^\frac{3}{4}}{2} \bigg\{ \io \deps \frac{w_{\eps x}^2}{g(\weps)} \bigg\}^\frac{1}{2}
	\|\psi\|_{L^2(\Omega)},
  \eea
  because $\frac{1}{\deps} \le \frac{1}{\sqrt{\eps}}$ in $\Omega$ according to (\ref{d_pos}).
  Furthermore,
  \be{32.4}
	\bigg|- \eps \io \sqrt{\deps} \frac{w_{\eps x}}{\sqrt{g(\weps)}} \psi_x  \bigg|
	\le \eps \bigg\{ \io \deps \frac{w_{\eps x}^2}{g(\weps)} \bigg\}^\frac{1}{2} \|\psi_x\|_{L^2(\Omega)},
  \ee
  whereas again invoking Lemma \ref{lem1} along with (\ref{g}) we see that
  \be{32.5}
	\bigg|- \io \sqrt{\deps} \frac{\ueps}{1+\eta_\eps u_\eps} g(\weps)\psi \bigg|
	\le g(M)\io \sqrt{\deps}u_\eps \psi 
	\le g(M) \bigg\{ \io \sqrt{\deps}^3 \ueps^3 \bigg\}^\frac{1}{3} \|\psi\|_{L^\frac{3}{2}(\Omega)}
  \ee
  for all $\eps\in (0,\eps_0)$. Thus, since $W^{1,2}(\Omega) \hra L^2(\Omega) \subset L^\frac{3}{2}(\Omega)$,
  and due to Lemma \ref{lem5} and Lemma \ref{lem11}, we obtain that for any $T>0$, 
  \bas
	\sup_{\eps\in (0,\es(T))} \int_0^T \bigg\{ \io \deps(x) \frac{w_{\eps x}^2(x,t)}{g(\weps(x,t))} dx 
		\bigg\}^\frac{3}{2} dt
	\le T \cdot \bigg\{ \sup_{\eps\in (0,\es(T))} \sup_{t\in (0,T)} 
	\io \deps(x) \frac{w_{\eps x}^2(x,t)}{g(\weps(x,t))} dx \bigg\}^\frac{3}{2} < \infty	
  \eas
  and
  \bas
	\sup_{\eps\in (0,\es(T))} \int_0^T \io \sqrt{\deps(x)}^3 \ueps^3(x,t) dxdt < \infty,
  \eas
  it follows from (\ref{32.3}), (\ref{32.4}), and (\ref{32.5}) that (\ref{32.2}) entails (\ref{32.1}).
\qed
\mysection{Construction of limit functions in $\{d>0\}$}\label{limits-d-positive}
We are now prepared for the construction of a limit function inside the positivity set of $d$ through 
a straightforward extraction process based on straighforward compactness arguments.
We remark that at this stage, besides the weighted functions $\sqrt{\deps}\weps$, our reasoning yet involves
the quantities $\sqrt{\deps(\ueps+1)}$, rather than those addressed in Lemma \ref{lem11}.
\begin{lem}\label{lem21}
  There exist a sequence $(\eps_k)_{k\in\N} \subset (0,\eps_0)$ and nonnegative functions 
  $\tu$ and $\tw$ defined in $\{d>0\}\times (0,\infty)$ such that $\eps_k\searrow 0$ as $k\to\infty$ and
  \begin{eqnarray}
	& & \ueps \to \tu
	\qquad \mbox{a.e.~in } \{d>0\} \times (0,\infty),
	\label{21.1} \\
	& & \ueps \wto \tu
	\qquad \mbox{in } L^1_{loc}([0,\infty);L^1(\{d>0\})),
	\label{21.2} \\
	& & \sqrt{\deps(\ueps+1)} \to \sqrt{d(\tu+1)}
	\qquad \mbox{in } L^2_{loc}([0,\infty);L^2(\{d>0\})),
	\label{21.3} \\
	& & \weps \to \tw
	\qquad \mbox{a.e.~in } \{d>0\} \times (0,\infty)
	\qquad \mbox{and} 
	\label{21.4} \\
	& & \sqrt{\deps} \weps \wto \sqrt{d} \tw
	\qquad \mbox{in } L^2_{loc}([0,\infty);W^{1,2}(\{d>0\}))
	\label{21.5}
  \end{eqnarray}
  as $\eps=\eps_k\searrow 0$.
\end{lem}
\proof
  Since given $T>0$ we can use (\ref{K1}) to estimate
  \bas
	\int_0^T \io \Big| (\sqrt{\deps(\ueps+1)})_x \Big|^2
	&=& \frac{1}{4} \int_0^T \io 
	\Big| \sqrt{\deps} \frac{u_{\eps x}}{\sqrt{\ueps+1}} + \frac{d_{\eps x}}{\sqrt{\deps}} \sqrt{\ueps+1} \Big|^2 \\
	&\le& \frac{1}{2} \int_0^T \io \deps \frac{u_{\eps x}^2}{\ueps+1}
	+ \frac{1}{2} \int_0^T \io \frac{d_{\eps x}^2}{\deps} (\ueps+1) \\
	&\le& \frac{1}{2} \int_0^T \io \deps \frac{u_{\eps x}^2}{\ueps+1}
	+ \frac{K_1}{2} \int_0^T \io (\ueps+1) 
  \eas
  for all $\eps\in (0,\eps_0)$, it follows from Lemma \ref{lem5} and Lemma \ref{lem2} that with $\es(T)$ as introduced
  there,
  \bas
	\Big(\sqrt{\deps(\ueps+1)}\Big)_{\eps\in (0,\es(T))}
	\quad \mbox{is bounded in } L^2((0,T);W^{1,2}(\Omega)).
  \eas
  Therefore, in view of Lemma \ref{lem31} the Aubin-Lions lemma asserts that for any such $T$,
  \bas
	\Big(\sqrt{\deps(\ueps+1)}\Big)_{\eps\in (0,\es(T))}
	\quad \mbox{is relatively compact in } L^2((0,T);L^2(\Omega)),
  \eas
  from which it follows by a standard argument that for a suitable sequence $(\eps_k)_{k\in\N} \subset (0,\eps_0)$
  and some $z\in L^2_{loc}([0,\infty);L^2(\Omega))$ we have $\eps_k\searrow 0$ as $k\to\infty$ and
  \be{21.6}
	\sqrt{\deps(\ueps+1)} \to z
	\qquad \mbox{in } L^2_{loc}([0,\infty);L^2(\Omega))
  \ee
  and
  \be{21.7}
	\sqrt{\deps(\ueps+1)} \to z
	\qquad \mbox{a.e.~in } \Omega\times (0,\infty)
  \ee
  as $\eps=\eps_k\searrow 0$.
  Since $\deps\to d$ a.e.~in $\Omega$ as $\eps\searrow 0$ by (\ref{d_conv}), this means that if we let 
  $\tu(x,t):=\frac{z^2(x,t)}{d(x)}-1$ for $x\in \{d>0\}$ and $t>0$, then (\ref{21.7}) and (\ref{21.6}) imply (\ref{21.1})
  and (\ref{21.3}), whereupon (\ref{21.1}) a posteriori also shows that $\tu$ must be nonnegative.\\
  We next make use of the estimate (\ref{5.1}) from Lemma \ref{lem5} to infer that for $T>0$ and $\es(T)$ as above,
  \bas
	\Big(\ueps \ln \ueps\Big)_{\eps\in (0,\es(T))}
	\quad \mbox{is bounded in } L^1(\Omega\times (0,T)),
  \eas
  so that the Dunford-Pettis theorem guarantees that $(\ueps)_{\eps\in (0,\es(T))}$ is relatively compact with respect
  to the weak topology in $L^1(\Omega\times (0,T))$, and that hence (\ref{21.2}) can be achieved on extracting a 
  subsequence of $(\eps_k)_{k\in\N}$ if necessary.\\
  As for the second solution component, we first use (\ref{K1}) to see that
  \bas
	\io \Big| (\sqrt{\deps} \weps)_x \Big|^2
	&=& \io \Big|\sqrt{\deps} w_{\eps x} + \frac{d_{\eps x}}{2\sqrt{\deps}} \weps \Big|^2 \\
	&\le& 2\io \deps w_{\eps x}^2 + \frac{1}{2} \io \frac{d_{\eps x}^2}{\deps} \weps^2 \\
	&\le& 2 \io \deps w_{\eps x}^2
	+ \frac{K_1}{2} \io \weps^2
	\qquad \mbox{for all } t>0,
  \eas
  so that for $T>0$ and $\es(T)$ as before, Corollary \ref{cor55} and Lemma \ref{lem1} warrant that
  \be{21.77}
	\Big(\sqrt{\deps} \weps\Big)_{\eps\in (0,\es(T))}
	\quad \mbox{is bounded in } L^\infty((0,T);W^{1,2}(\Omega)).
  \ee
  Thus,
  \be{21.8}
	\Big(\sqrt{\deps} \weps\Big)_{\eps\in (0,\es(T))}
	\quad \mbox{is relatively compact with respect to the weak topology in } L^2((0,T);W^{1,2}(\Omega)),
  \ee
  whereas (\ref{21.77}) in conjunction with Lemma \ref{lem32} and the Aubin-Lions lemma ensures that
  \bas
	\Big(\sqrt{\deps} \weps\Big)_{\eps\in (0,\es(T))}
	\quad \mbox{is relatively compact in } L^2(\Omega\times (0,T)).
  \eas
  Consequently, arguing as above we conclude upon passing to a further subsequence if necessary that both
  (\ref{21.4}) and (\ref{21.5}) hold with some $\tw: \{d>0\} \times (0,\infty)\to [0,\infty)$.
\qed
In dealing with the taxis term in (\ref{0eps}), the following consequence of (\ref{21.5}) will turn out to 
be more convenient.
\begin{cor}\label{cor211}
  With $(\eps_k)_{k\in\N}\subset (0,\eps_0)$ and $\tw$ as in Lemma \ref{lem21}, we have
  \be{211.1}
	\sqrt{\deps} w_{\eps x} \wto \sqrt{d} \tw_x
	\qquad \mbox{in } L^2_{loc}([0,\infty);L^2(\{d>0\}))
  \ee
  as $\eps=\eps_k\searrow 0$.
\end{cor}
\proof
  We rewrite
  \be{211.2}
	\sqrt{\deps} w_{\eps x}
	= (\sqrt{\deps} \weps)_x - \frac{d_{\eps x}}{2\sqrt{\deps}} \weps
  \ee
  and note that in view of the dominated convergence theorem, 
  combining (\ref{dx_conv}), (\ref{d_conv}) and (\ref{21.4}) with (\ref{K1}) and Lemma \ref{lem1} shows that
  for any $T>0$ and $\varphi\in L^2(\{d>0\} \times (0,T))$ we obtain
  \bas
	\int_0^T \io \int_{\{d>0\}} \frac{d_{\eps x}}{2\sqrt{\deps}} \weps \varphi
	\to \int_0^T \int_{\{d>0\}} \frac{d_x}{2\sqrt{d}} \tw \varphi
  \eas
  and hence
  \bas
	\frac{d_{\eps x}}{2\sqrt{\deps}} \weps \wto \frac{d_x}{2\sqrt{d}} \tw
	\qquad \mbox{in } L^2_{loc}([0,\infty);L^2(\{d>0\}))
  \eas
  as $\eps=\eps_k\searrow 0$.
  Therefore, (\ref{211.1}) results from (\ref{211.2}) on using (\ref{21.5}).
\qed
\mysection{Further convergence and integrability properties}\label{further-limits}
Let us now make use of the pointwise convergence property (\ref{5.1}) from Lemma \ref{lem5} to accomplish
our previously formulated goal concerning strong $L^2$ convergence of 
$\sqrt{\deps} \frac{\ueps}{(1+\eeps\ueps)^2}$.
Indeed, through an argument based on Egorov's theorem this will result from the fact that Lemma \ref{lem11}
implies bounds for this quantity in Lebesgue spaces involving superquadratic integrability.
\begin{lem}\label{lem12}
  Let $(\eps_k)_{k\in\N}\subset (0,\eps_0)$ be as provided by Lemma \ref{lem21}. Then
  \be{12.1}
	\sqrt{\deps} \frac{\ueps}{(1+\eeps\ueps)^2} \to \sqrt{d}\tu
	\qquad \mbox{in } L^2_{loc}([0,\infty);L^2(\{d>0\}))
  \ee
  as $\eps=\eps_k\searrow 0$.
\end{lem}
\proof
  As a consequence of Lemma \ref{lem11}, given $T>0$ we can find $c_1=c_1(T)>0$ such that with $\es(T)$ as in 
  Lemma \ref{lem5},
  \be{12.2}
	\int_0^T \io \sqrt{\deps}^3 \ueps^3 \le c_1
	\qquad \mbox{for all } \eps\in (0,\es(T)).
  \ee
  Since $\eeps>0$ for all $\eps\in (0,\eps_0)$, this implies that for 
  \bas
	\zeps:=\sqrt{\deps} \frac{\ueps}{(1+\eeps\ueps)^2},
	\qquad \eps\in (0,\es(T)),
  \eas
  we have
  \be{12.3}
	\int_0^T \io \zeps^3 \le c_1
	\qquad \mbox{for all } \eps\in (0,\es(T)).
  \ee
  Since $\eeps\to 0$ as $\eps\searrow 0$ by (\ref{eta_conv}), from Lemma \ref{lem21} we moreover know that
  \be{12.33}
	\zeps \to \sqrt{d}\tu
	\quad \mbox{a.e.~in } \{d>0\}\times (0,\infty)
	\qquad \mbox{as } \eps=\eps_k\searrow 0.
  \ee
  Therefore, according to a standard argument involving Egorov's theorem it particularly follows from (\ref{12.3}) 
  that
  \bas
	\zeps \wto \sqrt{d}\tu
	\quad \mbox{in } L^2(\{d>0\} \times (0,T))
	\qquad \mbox{as } \eps=\eps_k\searrow 0,
  \eas
  so that it remains to show that
  \be{12.4}
	\limsup_{\eps=\eps_k\searrow 0} \int_0^T \int_{\{d>0\}} \zeps^2 
	\le \int_0^T \int_{\{d>0\}} d\tu^2.
  \ee
  To this end, supposing on the contrary that for some $c_2 > \int_0^T \int_{\{d>0\}} du^2$ and some subsequence
  $(\eps_{k_j})_{j\in\N}$ of $(\eps_k)_{k\in\N}$ we had
  \be{12.5}
	\int_0^T \int_{\{d>0\}} \zeps^2 \to c_2
	\qquad \mbox{as } \eps=\eps_{k_j}\searrow 0,
  \ee
  once more by means of (\ref{12.3}) we could extract a further subsequence, again denoted by $(\eps_{k_j})_{j\in\N}$
  here fore convenience,
  along which for some $\widehat{z}\in L^\frac{3}{2}(\{d>0\} \times (0,T))$ we would have
  \bas
	\zeps^2 \wto \widehat{z}
	\quad \mbox{in } L^\frac{3}{2}(\{d>0\}\times (0,T))
	\qquad \mbox{as } \eps=\eps_{k_j}\searrow 0.
  \eas
  Since (\ref{12.33}) warrants that $\zeps^2 \to du^2$ a.e.~in $\{d>0\} \times (0,\infty)$ as $\eps=\eps_k\searrow 0$,
  again by Egorov's theorem this would imply that actually 
  \bas
	\zeps^2 \wto d\tu^2
	\quad \mbox{in } L^\frac{3}{2}(\{d>0\}\times (0,T))
	\qquad \mbox{as } \eps=\eps_{k_j}\searrow 0,
  \eas
  so that since the boundedness of $\{d>0\}\times (0,T)$ allows for choosing nontrivial constants as test functions here,
  we would conclude that
  we would conclude that
  \bas
	\int_0^T \int_{\{d>0\}} \zeps^2 \to \int_0^T \int_{\{d>0\}} du^2
	\qquad \mbox{as } \eps=\eps_{k_j} \searrow 0.
  \eas
  This contradiction to (\ref{12.5}) shows that in fact (\ref{12.4}) must hold, 
   whence the proof becomes complete.
\qed
A further property of the limit couple $(\tu,\tw)$, quite plausible in view of Corollary \ref{cor55}, can also be justified
on the basis of Egorov's theorem.
\begin{lem}\label{lem23}
  Suppose that $\tu$ and $\tw$ are as constructed in Lemma \ref{lem21}. Then for all $T>0$,
  \be{23.1}
	\int_0^T \int_{\{d>0\}} d\tu \tw_x^2 < \infty.
  \ee
\end{lem}
\proof
  According to Lemma \ref{lem21}, (\ref{eta_conv}), and Corollary \ref{cor211}, with $(\eps_k)_{k\in\N} \subset (0,\eps_0)$
  as in Lemma \ref{lem21} we have
  \be{23.2}
	\sqrt{\frac{\ueps}{1+\eeps\ueps}} \to \sqrt{\tu}
	\qquad \mbox{a.e.~in } \{d>0\} \times (0,T)
  \ee
  and
  \be{23.3}
	\sqrt{\deps} w_{\eps x} \wto \sqrt{d} \tw_x
	\qquad \mbox{in } L^2(\{d>0\}\times (0,T))
  \ee
  as $\eps=\eps_k\searrow 0$.
  Next, Corollary \ref{cor55} entails that with $\es(T)\in (0,\eps_0)$ taken from Lemma \ref{lem5}, the family
  $\Big(\sqrt{\deps\frac{\ueps}{1+\eeps\ueps}} w_{\eps x}\Big)_{\eps\in (0,\es(T))}$ is bounded
  in $L^2(\{d>0\} \times (0,T))$, so that we can find 
  \be{23.33}
	z\in L^2(\{d>0\}\times (0,T))
  \ee
  and a subsequence $(\eps_{k_j})_{j\in\N}$ of $(\eps_k)_{k\in\N}$ in such a way that
  \bas
	\sqrt{\frac{\ueps}{1+\eeps\ueps}} \cdot \Big(\sqrt{\deps}w_{\eps x}\Big)
	\equiv \sqrt{\deps\frac{\ueps}{1+\eeps\ueps}} w_{\eps x}
	\wto z
	\qquad \mbox{in } L^2(\{d>0\}\times (0,T))
  \eas
  as $\eps=\eps_{k_j} \searrow 0$. Here a known consequence of Egorov's theorem (\cite[Lemma A.1]{zhigun_surulescu_uatay})
  asserts that due to (\ref{23.2}) and (\ref{23.3}) we may identify
  \bas
	z=\sqrt{\tu} \cdot \Big(\sqrt{d} \tw_x\Big) \equiv \sqrt{d\tu} \tw_x
	\qquad \mbox{a.e.~in } \{d>0\}\times (0,T),
  \eas
  so that (\ref{23.1}) results from (\ref{23.33}).
\qed

\mysection{Solution properties of $\tu$ and $\tw$}\label{sol-prop}
We are now ready to make sure that $(\tu,\tw)$ indeed solves (\ref{0}) when restricted to $\{d>0\}$ in the following 
sense.
\begin{lem}\label{lem34}
  Let $\tu$ and $\tw$ be as obtained in Lemma \ref{lem21}.\abs
  i) \ If $\varphi \in C_0^\infty(\bar\Omega\times [0,\infty))$ is such that $\varphi_x=0$
  on $\pO\times (0,\infty)$ and additionally
  \be{34.01}
	\supp \varphi \subset \{d>0\} \times [0,\infty),
  \ee
  then 
  \bea{34.1}
	& & \hspace*{-20mm}
	- \int_0^\infty \int_{\{d>0\}} \tu \varphi_t
	- \int_{\{d>0\}} u_0\varphi(\cdot,0)
	= \int_0^\infty \int_{\{d>0\}} d\tu \varphi_{xx}
	+ \int_0^\infty \int_{\{d>0\}} d\tu \tw_x \varphi_x \nn\\
	& & + \int_0^\infty \int_{\{d>0\}} \tu f(\cdot,\tu,\tw) \varphi.
  \eea
  ii) \ For all $\varphi\in C_0^\infty(\Omega\times [0,\infty))$ fulfilling (\ref{34.01}), we have
  \be{34.11}
	\int_0^\infty \int_{\{d>0\}} \tw \varphi_t
	+ \int_{\{d>0\}} w_0 \varphi(\cdot,0)
	= \int_0^\infty \tu g(\tw) \varphi.
  \ee
\end{lem}
\proof
  On testing the first equation in (\ref{0eps}) by $\varphi$ we see that
  \bea{34.2}
	& & \hspace*{-20mm}
	- \int_0^\infty \io \ueps \varphi_t
	- \io u_0 \varphi(\cdot,0)
	= \int_0^\infty \io \deps \ueps \varphi_{xx}
	+\int_0^\infty \io \deps \frac{\ueps}{(1+\eeps\ueps)^2} w_{\eps x} \varphi_x \nn\\
	& & + \int_0^\infty \io \ueps f(\cdot,\ueps,\weps) \varphi
  	\qquad \mbox{for all } \eps \in (0,\eps_0),
  \eea
  where since $\ueps\wto \tu$ in $L^1_{loc}([0,\infty);L^1(\{d>0\}))$ as $\eps=\eps_k\searrow 0$ by Lemma \ref{lem21},
  according to (\ref{34.01}) we have
  \bas
	- \int_0^\infty \io \ueps\varphi_t \to - \int_0^\infty \int_{\{d>0\}} \tu \varphi_t
  \eas
  and
  \bas
	\int_0^\infty \io \deps\ueps \varphi_{xx}
	\to \int_0^\infty \int_{\{d>0\}} d\tu \varphi_{xx}
  \eas
  as $\eps=\eps_k\searrow 0$, because $\deps\to d$ in $L^\infty(\Omega)$ as $\eps\searrow 0$ due to (\ref{d_conv}).\\
  Next, since Lemma \ref{lem21} warrants that also $\ueps \to \tu$ and $\weps \to \tw$ a.e.~in $\{d>0\} \times (0,\infty)$
  as $\eps=\eps_k\searrow 0$, it follows from Lemma \ref{lem6} and a standard argument,
  again involving Egorov's theorem, that
  \bas
	\ueps f(\cdot,\ueps,\weps) \wto \tu f(\cdot,\tu,\tw)
	\qquad \mbox{in } L^1_{loc}([0,\infty);L^1(\{d>0\})),
  \eas
  and that hence
  \bas
	\int_0^\infty \io \ueps f(x,\ueps,\weps) \to \int_0^\infty \int_{\{d>0\}} \tu f(x,\tu,\tw)
  \eas
  as $\eps=\eps_k\searrow 0$.\\
  Finally, from Corollary \ref{cor211} we know that
  \bas
	\sqrt{\deps} w_{\eps x} \wto \sqrt{d} \tw_x
	\qquad \mbox{in } L^2_{loc}([0,\infty);L^2(\{d>0\}))
  \eas
  as $\eps=\eps_k\searrow 0$, which combined with the strong convergence property of
  $\sqrt{\deps} \frac{\ueps}{(1+\eeps\ueps)^2}$ in $L^2_{loc}([0,\infty);L^2(\{d>0\}))$ asserted by Lemma \ref{lem12}
  ensures that
  \bas
	\int_0^\infty \io \deps \frac{\ueps}{(1+\eeps\ueps)^2} w_{\eps x} \varphi_x
	= \int_0^\infty \io \Big(\sqrt{\deps} \frac{\ueps}{(1+\eeps\ueps)^2}\Big)\cdot\Big(\sqrt{\deps} w_{\eps x}\Big) 
		\varphi_x 
	\to \int_0^\infty \int_{\{d>0\}} d\tu \tw_x \varphi_x
  \eas
  as $\eps=\eps_k\searrow 0$. Therefore, (\ref{34.1}) is a consequence of (\ref{34.2}).\abs
  To verify (\ref{34.11}), given $\varphi\in C_0^\infty(\Omega\times [0,\infty))$ fulfilling (\ref{34.01}) we obtain from
  (\ref{0eps}) that
  \be{34.3}
	\int_0^\infty \io \weps \varphi_t + \io w_{0\eps} \varphi(\cdot,0)
	= - \eps \int_0^\infty \io \frac{w_{\eps x}}{\sqrt{g(\weps)}} \varphi_x
	- \int_0^\infty \io \frac{\ueps}{1+\eta _\eps\ueps} g(\weps) \varphi 
  \ee
  for all $\eps\in (0,\eps_0)$. Here by Lemma \ref{lem21}, Lemma \ref{lem1} and the dominated convergence theorem,
  \be{34.4}
	\int_0^\infty \io \weps \varphi_t \to \int_0^\infty \int_{\{d>0\}} \tw \varphi_t
  \ee
  as $\eps=\eps_k\searrow 0$, whereas (\ref{w0eps}) trivially ensures that
  \be{34.5}
	\io w_{0\eps} \varphi(\cdot,0)
	\to \int_{\{d>0\}} w_0 \varphi(\cdot,0)
  \ee
  as $\eps=\eps_k\searrow 0$.
  Moreover, combining Lemma \ref{lem6} with the pointwise convergence properties in (\ref{21.1}) and (\ref{21.4}) we
  easily infer that
  \bas
	\ueps g(\weps) \wto \tu g(\tw)
	\qquad \mbox{in } L^1_{loc}([0,\infty);L^1(\{d>0\}))
  \eas
  and thus, also relying on (\ref{eta_conv}) and again on Lemma \ref{lem21}, we obtain
  \be{34.6}
	\int_0^\infty \io \frac{\ueps}{1+\eta_\eps \ueps} g(\weps)\varphi
	\to \int_0^\infty \int_{\{d>0\}} \tu g(\tw)\varphi
  \ee
  as $\eps=\eps_k\searrow 0$.
  Finally, once more relying on the fact that $\deps\ge \sqrt{\eps}$ by (\ref{d_pos}), we see by using the Cauchy-Schwarz
  inequality that if $T>0$ is large enough such that $\varphi\equiv 0$ in $\Omega\times (T,\infty)$, then
  \bas
	\bigg| \eps \int_0^\infty \io \frac{w_{\eps x}}{\sqrt{g(\weps)}} \varphi_x \bigg|
	&\le& \eps \int_0^\infty \bigg\{ \io \deps \frac{w_{\eps x}^2}{g(\weps)} \bigg\}^\frac{1}{2} \cdot
	\bigg\{ \io \frac{\varphi_x^2(\cdot,t)}{\deps} \bigg\}^\frac{1}{2} dt \\
	&\le& \eps^\frac{3}{4} \sup_{t\in (0,T)} \bigg\{ \io \deps \frac{w_{\eps x}^2(\cdot,t)}{g(\weps(\cdot,t))} 
		\bigg\}^\frac{1}{2}
	\cdot \int_0^\infty \bigg\{ \io \varphi_x^2(\cdot,t) \bigg\}^\frac{1}{2} dt
  \eas
  for all $\eps\in (0,\eps_0)$. 
  Therefore, since with $\es(T)\in (0,\eps_0)$ as given by Lemma \ref{lem5} we know from (\ref{5.2}) that
  \bas
	\sup_{\eps\in (0,\es(T))} \sup_{t\in (0,T)} \io \deps \frac{w_{\eps x}^2(\cdot,t)}{g(\weps(\cdot,t))} < \infty,
  \eas
  it follows that
  \bas
	\eps\int_0^\infty \io \frac{w_{\eps x}}{\sqrt{g(\weps)}} \varphi_x \to 0
  \eas
  as $\eps=\eps_k\searrow 0$.
  In combination with (\ref{34.4})-(\ref{34.6}), this shows that (\ref{34.3}) implies (\ref{34.11}).
\qed
\mysection{Proof of Theorem \ref{theo33}}\label{the-proof}
We can finally extend the above spatially local solution in an evident manner so as to become
a global weak solution in the flavor of Definition \ref{defi_weak}. 
In the verification of the desired solution property near the boundary of $\{d>0\}$ we shall make use of the following
consequence of the inclusion $\sqrt{d} \in W^{1,\infty}(\Omega)$.
\begin{lem}\label{lem35}
  Let $x\in\bar\Omega$. Then
  \be{35.1}
	d(x) \le \frac{K_1}{4} \Big\{ \dist (x, \{d=0\}) \Big\}^2.
  \ee
\end{lem}
\proof
  We only need to consider the case when $d(x)>0$, in which by the closedness of $\{d=0\}$ we can pick $x_0\in\bar\Omega$
  such that $d(x_0)=0$ and $|x-x_0|=\dist (x,\{d=0\})>0$. 
  Thanks to (\ref{K1}), we then have
  \bas
	\sqrt{d(x)} = \int_{x_0}^x (\sqrt{d})_y dy
	= \frac{1}{2} \int_{x_0}^x \frac{d_x(y)}{\sqrt{d(y)}} dy
	\le \frac{\sqrt{K_1}}{2} |x-x_0|
	= \frac{\sqrt{K_1}}{2} \cdot \dist (x,\{d=0\}),
  \eas
  from which (\ref{35.1}) follows.
\qed
By means of an appropriate cut-off procedure we can thereby proceed to show that the natural extension
of $(\tu,\tw)$, consisting of a solution to the ODE system formally associated with (\ref{0}) in $\{d=0\}$ 
indeed solves (\ref{0}) in the desired sense.\abs
\proofc of Theorem \ref{theo33}. \quad
  We let $\tu$ and $\tw$ denote the functions defined on $\{d>0\} \times (0,\infty)$ in Lemma \ref{lem21}, and for
  fixed $x\in \{d=0\}$ we let $(\hu(x,\cdot),\hw(x,\cdot)) \in (C^1([0,\infty)))^2$ be the solution of the 
  initial-value problem 
  \be{33.2}
	\left\{ \begin{array}{l}
	\hu_t=\hu f(x,\hu,\hw), \qquad t>0, \\[1mm]
	\hw_t=-\hu g(\hw), \qquad t>0, \\[1mm]
	\hu(x,0)=u_0(x), \quad \hw(x,0)=w_0(x).
	\end{array} \right.
  \ee
  Indeed, it follows from (\ref{dfg_reg}), (\ref{f_upper}) and (\ref{g}) that for any such $x$ this ODE problem
  possesses a globally defined solution fulfilling
  \be{33.3}
	0 \le \hw(x,t) \le M
	\qquad \mbox{for all } t>0
  \ee
  and
  \be{33.4}
	0 \le \hu(x,t) \le u_0(x) e^{\rho(M)t}
	\qquad \mbox{for all } t>0,
  \ee
  and since $u_0$ and $w_0$ are continuous in $\bar\Omega$ by (\ref{init}), standard ODE theory warrants that
  both $\hu$ and $\hw$ are continuous in $\{d=0\} \times [0,\infty)$. Therefore,
  \be{33.44}
	(u,w)(x,t):= \left\{ \begin{array}{ll}
	(\tu,\tw)(x,t), \qquad & x\in \{d>0\}, \ t>0, \\[1mm]
	(\hu,\hw)(x,t), \qquad & x\in \{d=0\}, \ t>0, 
	\end{array} \right.
  \ee
  defines a pair of nonnegative measurable functions on all of $\Omega\times (0,\infty)$ which thanks to Lemma \ref{lem21},
  Lemma \ref{lem6}, Lemma \ref{lem1}, (\ref{g}),
  (\ref{33.3}) and (\ref{33.4}) satisfy (\ref{w1}) and (\ref{w11}), and for which Lemma \ref{lem23} in particular entails 
  that also (\ref{w2}) holds.\abs
  In order to verify (\ref{w3}), we first make use of the fact that by continuity of $d$ 
  the set $\{d>0\}$ is relatively open in $\bar\Omega$, 
  and hence consists of countably many connected components; that is, there exist an index set
  $I\subset\N$ and intervals $P_i\subset\bar\Omega$, $i\in I$, such that $\{d>0\}=\bigcup_{i\in I} P_i$
  and $P_i \cap P_j =\emptyset$ for $i,j\in I$ with $i\ne j$. 
  Now for each $i\in I$, there exist $a_i\in\bar\Omega$ and $b_i\in\bar\Omega$ such that
  $(a_i,b_i)\subset P_i \subset [a_i,b_i]$, where $a_i\in P_i$ if and only if $a_i\in \pO$ and $b_i\in P_i$
  if and only if $b_i\in\pO$. 
  For fixed $\delta\in (0,1)$, defining $\delta_i:=2^{-i}\delta$, $i\in I$, it is then possible to pick
  $(\zd^{(i)})_{i\in I} \subset C^\infty(\bar\Omega)$ such that for all $i\in I$ we have $0\le \zd^{(i)} \le 1$ in
  $\bar\Omega$, $\zd^{(i)}(x)=1$ whenever $x\in P_i$ is such that $\dist(x,\partial P_i)\ge \delta_i$,
  $\zd^{(i)}\equiv 0$ in $\bar\Omega\setminus P_i$, and
  \be{33.55}
	|\zeta_{\delta x}^{(i)}| \le \frac{2}{\delta_i}
	\qquad \mbox{in } \bar\Omega
  \ee
  as well as
  \be{33.56}
	|\zeta_{\delta xx}^{(i)}| \le \frac{16}{\delta_i^2}
	\qquad \mbox{in } \bar\Omega,
  \ee
  where in the exceptional case $a_i\in\pO$ we can additionally achieve that $\zd^{(i)} \equiv 1$ holds even throughout
  $[a_i,b_i-\delta_i]$, and where, similarly, in the case $b_i\in\pO$ we require that $\zd^{(i)} \equiv 1$
  in $[a_i+\delta_i,b_i]$.\\
  Now given $\varphi\in C_0^\infty(\bar\Omega\times [0,\infty))$ satisfying $\varphi_x=0$ on $\pO\times (0,\infty)$,
  from (\ref{33.2}) and (\ref{33.4}) we immediately see that
  \be{33.6}
	- \int_0^\infty \int_{\{d=0\}} u \varphi_t - \int_{\{d=0\}} u_0 \varphi(\cdot,0)
	= \int_0^\infty \int_{\{d=0\}} uf(x,u,w) \varphi.
  \ee
  Moreover, Lemma \ref{lem34} guarantees that if we let
  \bas
	\zd:=\sum_{i\in I} \zd^{(i)}, 
	\qquad \delta\in (0,1),
  \eas
  then since $\supp (\zd \cdot \varphi) \subset \{d>0\} \times [0,\infty)$, we have
  \bea{33.7}
	& & \hspace*{-20mm}
	- \int_0^\infty \int_{\{d>0\}} \zd u\varphi_t
	-\int_{\{d>0\}} \zd u_0 \varphi(\cdot,0) \nn\\
	&=& \int_0^\infty \int_{\{d>0\}} du \cdot (\zd \varphi)_{xx}
	+ \int_0^\infty \int_{\{d>0\}} duw_x \cdot (\zd \varphi)_x \nn\\
	& & + \int_0^\infty \int_{\{d>0\}} \zd u f(\cdot,u,w)\varphi \nn\\
	&=& \int_0^\infty \int_{\{d>0\}} \zd du\varphi_{xx} 
	+ 2\int_0^\infty \int_{\{d>0\}} \zeta_{\delta x} du \varphi_x
	+ \int_0^\infty \int_{\{d>0\}} \zeta_{\delta xx} du\varphi \nn\\
	& & + \int_0^\infty \int_{\{d>0\}} \zd duw_x \varphi_x
	+ \int_0^\infty \int_{\{d>0\}} \zeta_{\delta x} duw_x \varphi \nn\\
	& & + \int_0^\infty \int_{\{d>0\}} \zd uf(\cdot,u,w)\varphi
	\qquad \mbox{for all } \delta\in (0,1).
  \eea
  Here we may use that $0\le\zd\le 1$ and that as $\delta\searrow 0$ we have $\zd \to 1$ a.e.~in $\{d>0\}$ to infer from
  the dominated convergence theorem that
  \be{33.8}
	- \int_0^\infty \int_{\{d>0\}} \zd u\varphi_t
	\to - \int_0^\infty \int_{\{d>0\}} u\varphi_t
  \ee
  and
  \be{33.9}
	- \int_{\{d>0\}} \zd u_0 \varphi(\cdot,0)
	\to - \int_{\{d>0\}} u_0 \varphi(\cdot,0)
  \ee
  as well as
  \be{33.10}
	\int_0^\infty \int_{\{d>0\}} \zd du \varphi_{xx}
	\to \int_0^\infty \int_{\{d>0\}} du\varphi_{xx}
  \ee
  and
  \be{33.11}
	\int_0^\infty \int_{\{d>0\}} \zd duw_x \varphi_x
	\to \int_0^\infty \int_{\{d>0\}} duw_x\varphi_x
  \ee
  and
  \be{33.12}
	\int_0^\infty \int_{\{d>0\}} \zd uf(\cdot,u,w)\varphi
	\to \int_0^\infty \int_{\{d>0\}} uf(\cdot,u,w)\varphi
  \ee
  as $\delta\searrow 0$.
  In order to estimate the integrals on the right of (\ref{33.7}) which contain derivatives of $\zd$, let us first
  observe that as a consequence of (\ref{33.55}), (\ref{33.56}) and Lemma \ref{lem35} we know that whenever
  $x\in\bar\Omega$ is such that $\zeta_{\delta x}(x)\ne 0$, for some $i\in I$ we have $x\in P_i$ and 
  $\dist(x,\{d=0\})\le \delta_i$ and hence
  \be{33.13}
	d(x)\zeta_{\delta x}^2(x) \le \frac{K_1 \delta_i^2}{4} \cdot \Big(\frac{2}{\delta_i}\Big)^2 = K_1
  \ee
  as well as
  \be{33.14}
	d(x) \cdot |\zeta_{\delta xx}(x)| \le \frac{K_1 \delta_i^2}{4} \cdot \frac{16}{\delta_i^2} = 4K_1.
  \ee
  Furthermore, again by mutual disjointness of the $P_i$,
  \bas
	\Big| \supp \zeta_{\delta x} \Big|
	\le \sum_{i\in I} 2\cdot\delta_i
	= \sum_{i\in I} 2\cdot (2^{-i} \delta) 
	\le 2\delta
	\qquad \mbox{for all } \delta\in (0,1),
  \eas
  so that since we know from Lemma \ref{lem11}, Lemma \ref{lem23} and Fatou's lemma that with $T>0$ taken large enough 
  fulfilling $\varphi\equiv 0$ in $\Omega\times (T,\infty)$ we have
  \bas
	\sqrt{d}^3 u^3 \in L^1(\{d>0\}\times (0,T))
  \eas
  and
  \bas
	duw_x^2 \in L^1(\{d>0\}\times (0,T)), 
  \eas
  from the dominated convergence theorem it follows that
  \be{33.15}
	\int_0^T \int_{\supp \zeta_{\delta x}} \sqrt{d}^3 u^3 \to 0
  \ee
  and
  \be{33.16}
	\int_0^T \int_{\supp \zeta_{\delta x}} duw_x^2 \to 0
  \ee
  as $\delta\searrow 0$, whereas combining (\ref{33.14}) with the dominated convergence theorem shows that also
  \be{33.17}
	\int_{\supp \zeta_{\delta x}} |\zeta_{\delta xx}| \cdot d \to 0
  \ee
  as $\delta\searrow 0$. 
 
  Thus, using the H\"older inequality along with (\ref{33.13}) and (\ref{33.15}) we obtain that
  \bea{33.18}
	\bigg|2 \int_0^\infty \int_{\{d>0\}} \zeta_{\delta x} du\varphi_x \bigg|
	&\le& 2\|\varphi_x\|_{L^\infty(\Omega\times (0,\infty))} 
	\bigg\{ \int_0^T \int_{\supp \zeta_{\delta x}} \sqrt{d}^3 u^3 \bigg\}^\frac{1}{3}
	\cdot \bigg\{ \int_0^T \int_{\supp \zeta_{\delta x}} |\zeta_{\delta x}|^\frac{3}{2} 
		d^\frac{3}{4} \bigg\}^\frac{2}{3} \nn\\
	&\le& 2\sqrt{K_1} T^\frac{2}{3} \|\varphi_x\|_{L^\infty(\Omega\times (0,\infty))}
	\bigg\{ \int_0^T \int_{\supp \zeta_{\delta x}} \sqrt{d}^3 u^3 \bigg\}^\frac{1}{3} \nn\\[2mm]
	&\to& 0
  \eea
  as $\delta\searrow 0$, while from (\ref{33.17}) we infer that
  \bea{33.19}
	\bigg| 	\int_0^\infty \int_{\{d>0\}} \zeta_{\delta xx} du\varphi \bigg|
	&\le& \|\varphi\|_{L^\infty(\Omega\times (0,\infty))}
	\int_0^T \Big\|\sqrt{d} u(\cdot,t)\Big\|_{L^\infty(\{d>0\})}
	\int_{\supp \zeta_{\delta x}} |\zeta_{\delta xx}| \cdot d \ dt \nn\\
	&\le& \sqrt{T} \|\varphi\|_{L^\infty(\Omega\times (0,\infty))}
	\bigg\{ \int_0^T \Big\|\sqrt{d} u(\cdot,t)\Big\|_{L^\infty(\{d>0\})} dt \bigg\}^\frac{1}{2} 
	\cdot \int_{\supp \zeta_{\delta x}} |\zeta_{\delta xx}| \cdot d \nn\\[2mm]
	&\to& 0
  \eea
  as $\delta\searrow 0$, because Lemma \ref{lem11} together with Fatou's lemma warrants that
  \bas
	\int_0^T \Big\|\sqrt{d}u(\cdot,t)\Big\|_{L^\infty(\{d>0\})}^2 dt
	\le \liminf_{\eps=\eps_k\searrow 0} \int_0^T \Big\|\sqrt{\deps}\ueps(\cdot,t)\Big\|_{L^\infty(\Omega)}^2 dt
	<\infty.
  \eas
  Since finally (\ref{33.16}) along with (\ref{33.13}) ensures that also
  \bas
	\bigg| \int_0^\infty \int{\{d>0\}} \zeta_{\delta x} duw_x \varphi \bigg|
	&\le& \|\varphi\|_{L^\infty(\Omega\times (0,\infty))}
	\bigg\{ \int_0^T \int_{\supp \zeta_{\delta x}} duw_x^2 \bigg\}^\frac{1}{2}
	\cdot \bigg\{ \int_0^T \io \zeta_{\delta x}^2 d \bigg\}^\frac{1}{2} \\
	&\le& \sqrt{K_1 T} \|\varphi\|_{L^\infty(\Omega\times (0,\infty))}
	\bigg\{ \int_0^T \int_{\supp \zeta_{\delta x}} duw_x^2 \bigg\}^\frac{1}{2} \nn\\[2mm]
	&\to& 0
  \eas
  as $\delta\searrow 0$, from (\ref{33.7})-(\ref{33.12}), (\ref{33.18}) and (\ref{33.19}) we conclude that
  \bas
	& & \hspace*{-20mm}
	- \int_0^\infty \int_{\{d>0\}} u\varphi_t
	- \int_{\{d>0\}} u_0 \varphi(\cdot,0)
	= \int_0^\infty \int_{\{d>0\}} du\varphi_{xx}
	+ \int_0^\infty \int_{\{d>0\}} duw_x \varphi_x \nn\\
	& & + \int_0^\infty \int_{\{d>0\}} uf(\cdot,u,w) \varphi,
  \eas
  which in combination with (\ref{33.6}) shows that indeed (\ref{w3}) is valid for any such $\varphi$.\abs
  The derivation of (\ref{w4}) is much less involved: Given $\varphi\in C_0^\infty(\Omega\times [0,\infty))$,
  from (\ref{33.44}) and (\ref{33.2}) we first obtain that
  \be{33.20}
	\int_0^\infty \int_{\{d=0\}} w\varphi_t
	+ \int_{\{d=0\}} w_0\varphi(\cdot,0)
	= \int_0^\infty \int_{\{d=0\}} ug(w)\varphi,
  \ee
  whereas with $(\zd)_{\delta\in (0,1)}$ as introduced above we obtain from Lemma \ref{lem34} that
  \be{33.21}
	\int_0^\infty \io \zd w\varphi_t
	+ \io \zd w_0 \varphi(\cdot,0)
	= \int_0^\infty \io \zd ug(w)\varphi
  \ee
  for all $\delta\in (0,1)$.
  Using that $w$ and $ug(w)$ belong to $L^1_{loc}([0,\infty);L^1(\{d>0\}))$ by Lemma \ref{lem1}, Lemma \ref{lem6}
  and Lemma \ref{lem21}, we may again employ the dominated convergence theorem here to see that in the limit 
  $\delta\searrow 0$, (\ref{33.21}) implies that
  \bas
	\int_0^\infty \int_{\{d>0\}} w\varphi_t
	+ \int_{\{d>0\}} w_0\varphi(\cdot,0)
	= \int_0^\infty \int_{\{d>0\}} ug(w)\varphi,
  \eas
  and that thus in view of (\ref{33.20}) also (\ref{w4}) holds.
\qed

\end{document}